\documentclass{amsart}

\usepackage[T1]{fontenc}
\usepackage[latin9]{inputenc}
\usepackage[margin=1in]{geometry}
\usepackage{amsthm}
\usepackage{amsmath}
\usepackage{amssymb}
\usepackage{amsrefs}
\usepackage{appendix}
\usepackage{color}
\usepackage{graphicx}
\usepackage{esint}
\usepackage{hyperref}
\usepackage{verbatim}

\usepackage{fourier}
\usepackage{xcolor}

\def\C{{\mathbb C}}
\def\R{{\mathbb R}}
\def\N{{\mathbb N}}
\def\Z{{\mathbb Z}}

\def\le{\leqslant}

\newcommand{\eps}{\varepsilon}
\newcommand{\e}{\varepsilon}
\newcommand{\wt}{\widetilde}

\theoremstyle{plain}
\newtheorem{theorem}{Theorem}[section]

\newtheorem*{hyp}{Assumption}

\theoremstyle{definition}

\newtheorem*{remark*}{Remark}

\numberwithin{equation}{section}

\allowdisplaybreaks

\title{Bloch dynamics with the second order Berry phase correction}

\author{Jianfeng Lu} 
\address{Department of Mathematics, Department of
  Physics, and Department of Chemistry, Duke University, Box 90320,
  Durham NC 27708, USA} 
\email{jianfeng@math.duke.edu}

\author{Zihang Zhang}
\address{School of Mathematical Sciences, Peking University, Beijing, 100871, China}
\email{zhang-zihang@pku.edu.cn}

\author{Zhennan Zhou}
\address{Beijing International Center for Mathematical Research, Peking University, Beijing, 100871, China}
\email{zhennan@bicmr.pku.edu.cn}

\date{First version: Dec 23, 2015; this version: \today}\thanks{This work is partially supported by the National
  Science Foundation under grants DMS-1312659, DMS-1454939 and
  RNMS11-07444 (KI-Net). J.L. is also partially supported by the
  Alfred P.~Sloan Foundation.  Z. Zhou is supported by NSFC grant No. 11801016. We thank Wolfgang Gaim and Stefan
  Teufel for helpful discussions.}

\begin{document}

\begin{abstract} { We derive the semiclassical Bloch
    dynamics with the second-order Berry phase correction in the
    presence of the slow-varying scalar potential as perturbation.
    Our mathematical derivation is based on a two-scale WKB asymptotic
    analysis.  For a uniform external electric field, the
    bi-characteristics system after a positional shift introduced by
    Berry connections agrees with the recent result in previous works.
    Moreover, for the case with a linear external electric field, we
    show that the extra terms arising in the bi-characteristics system
    after the positional shift are also gauge independent. }
\end{abstract}

\maketitle
\section{Introduction}

The understanding of dynamics of Bloch electrons and their response to
external electromagnetic fields plays an important role in solid state
physics (see for example \cites{BF, Berry, RS, XCN} and the references
therein). In recent years, many works such as \cites{Berry,
	Bloch, GYN, XCN, YuLiuYang, GaoXiao, SunLu, XiaoDuNiu, SunLu2} have explored the significant role of the Berry
phase in Bloch dynamics and vast related fields. 
There have been
series of important mathematical works in the direction as well, which
are devoted to rigorously justify the the validity of the physics
models and provide insight for possible generalizations (see for
example \cites{CMS, ELY, MMP, PST} and the references therein).

Under the single-particle approximation, the dynamics of an electron
is treated as an independent particle on the effective periodic
potential generated by ions and other electrons (as a mean-field) in
the crystal. After non-dimensionalization, the dynamics is given
by
\begin{equation}\label{SE1}
i\eps \frac{ \partial}{\partial t}\psi(t,x)=H \psi(t,x)= \left(-\frac{\eps^2}{2}\Delta_x+V\left( \frac{x}{\eps}\right)+U(x)\right)\psi(t,x),
\end{equation}
where $\psi: [0, \infty) \times \R^d \to \mathbb{C}$ is a single
particle wave function, $\eps$ is the semi-classical parameter, $V(z)$
is the lattice potential which is periodic with respect to $\mathbb{L}$,
and $U(x)$ is the slowly-varying scalar potential. 

{The dimensionless equation \eqref{SE1} is obtained from the original Schr\"odinger equation in physical units
\begin{equation}\label{SEphy}
	i\hbar\frac{\partial}{\partial t}\psi=-\frac{\hbar^2}{2m}\Delta \psi +V(x)\psi-U(x)\psi,
\end{equation}
following the procedures as in, e.g., \cite{ELY}. Here, $m$ is the mass and $\hbar$ is the reduced Planck constant. We introduce $l$ as the lattice constant, and $L$ and $T$ as the macroscopic length and time
scales, respectively. Following the calculations in Appendix \ref{app:a}), we obtain two dimensionless parameters $\e=l/L,\ h=\hbar T/(mL^2)$. In this paper we only consider the
distinguished limit when $\e =h\ll 1$. In such setting, the typical wavelength of the wave function is comparable with the size of the lattice, which are much smaller than the observation scales in space and time. The external potential $U(x)$ is slowly-varying because $U(\cdot)$ only depends on the macroscopic length scale $x$, not the quantum length scale $z=x/\e$. Such scaling has been vastly adopted in mathematical analysis of Bloch electrons, and the role of the slowly varying potential has been explored, see \cites{CMS, GMMP, SMM}.

}

It is well known that such semi-classical Schr\"odinger equations
propagate oscillations of order $\mathcal{O}(\e)$ both in space and
time.  With this model, the relevant physical scale translates to the
case when the typical wavelength is comparable to the period of the
medium, and both of which are assumed to be small on the length-scale
of the considered physical domain. This consequently leads us to a
problem involving two-scales where from now on we shall denote by $0 <
\e \ll 1$ the small dimensionless parameter describing the
microscopic/macroscopic scale ratio. We remark that, equation
\eqref{SE1} can also be derived from the Schr\"odinger equation in
physical units by introducing certain rescaling, which we shall omit in
this paper. The readers may refer to \cites{CMS, ELY} for such
calculations.


The electronic dynamics in crystals have been studied for many years
in the semi-classical regime, where the Liouville equations replace the role
of the Schr\"odinger equation in the limit when the rescaled Planck constant
tends to zero. With the help of the Bloch-Floquet theory \cite{RS}, Markowich,
Mauser and Poupaud in \cite{MMP} derived the semi-classical Liouville equation for
describing the propagation of the phase-space density for an energy band,
which controls the macroscopic dynamic behavior of the electrons. Later
these results were generalized to the case when a weak random potential in \cite{BF}
and in \cite{GMMP} nonlinear interactions were present.

Berry phase is an important object that appears during the adiabatic
limit of quantum dynamics, as some slow-changing variables enter the
quantum evolution as parameters, see \cites{Berry, SW}. As observed by
Simon in \cite{Simon}, the adiabatic Berry phase has an elegant
mathematical interpretation as the holonomy of a certain connection,
the Berry connection, in the appropriate fiber bundle. This setup
gives rise to the Berry curvature, which is gauge invariant and can be
considered as a physical observable. It has been used in the Bloch
dynamics to explain various important phenomena in crystals, see for
example \cites{GYN,XCN} and related references. Panati, Spohn and
Teufel later gave a rigorous derivation of such Bloch dynamics in
\cites{PST,PST2} by writing down the effective Hamiltonian with the
help of the Weyl quantization. A simple derivation of the Bloch
dynamics with Berry phase correction based on WKB asymptotics was also
given by E, Lu and Yang in \cite{ELY}.


Recently, Gao, Yang and Niu in
\cites{GYN,GYN2} constructed a second order semi-classical theory
for Bloch electrons under uniform electromagnetic fields, {which was based on the semiclassical Lagrangian  approach originated in \cite{SunNiu}}. The second
order correction terms in Bloch dynamics are obtained, where the first
order correction to the Berry curvature is derived. The second order
semi-classical theory can be used to explain some important physical
phenomenon, such as electric polarizability, magnetic susceptibility,
and magnetoelectric polarizability, etc., { which are investigated in the recent literature  \cites{YuLiuYang, GaoXiao, SunLu, XiaoDuNiu, SunLu2}.  In particular, some of the applications are connected with the second order Bloch dynamics with only the external electric field as perturbation, see e.g. \cites{GaoXiao, XiaoDuNiu}.} { This provides the motivation
to our current study: we aim to extend the mathematical derivation
of Bloch dynamics to include second order corrections  in the presence of an external electric field}.

The main purpose of this paper is to give a {rigorous} derivation of the
effective Bloch dynamics in crystals based on asymptotic analysis up to
the second
order. 
{The justification of the corrected dynamics is not based on the semiclassical Lagrangian  approach, but the WKB analysis. Since with the proper incorporation of the Bloch function, the WKB solution with high order corrections is proved in \cite{CMS} to approximate the true solution of \eqref{SE1} up to a quantifiable error in an arbitrary order of $\eps$,  this approximation theory acts as the foundation of our derivation.}
The final results that we obtain agree with the recent paper
\cite{GYN} in the situation of uniform electric field, but our
approach is mathematically rigorous and is able to handle more general
potentials.   We also
note an independent derivation of the Bloch dynamics with second order
correction \cite{GaimTeufel} using Weyl quantization of operator
valued symbols (see \cite{StiepanTeufel} for a related work).

More specifically, with a two-scaled WKB ansatz for equation
\eqref{SE1}, we derive the phase equation with the second order
corrections, and correspondingly the perturbation in Hamiltonian and
in Bloch energy.  The truncated WKB solution is proved to be a valid
approximation to the exact solution. {  The phase equation
  with the second order corrections can still be viewed as a
  Hamilton-Jacobi equation when the leading order phase equation and
  the amplitude equation are solved in the pre-processing stage. }
Note that, since we study here WKB type solutions to the Schr\"odinger
equation, the asymptotic solution we derive is valid only before
caustics. If long time validity of the asymptotic solution is desired,
one needs to consider instead for example the Gaussian beam methods
\cites{DGR, JWYH}, the Wigner functions \cites{LP,SMM}, or the frozen
Gaussian approximation for periodic media \cites{DeLuYa1,
  DeLuYa2}. The derivations of Bloch dynamics with Berry phase
corrections using these approaches are interesting future
directions. { We further derive the characteristic
  equations of the phase equation with the second order
  corrections. In two specific cases: when the electric field is
  uniform or linear in space, we perform a physics-inspired change of
  variable and show that the all the terms in the resulting equations
  are gauge independent except those from the extra wave packet
  energy. In the former case, our result essentially agrees with the
  recent results  \cites{GYN,GYN2} and thus
  provides a rigorous justification of these results. In the latter
  case of linear electric field, we find some additional gauge
  independent terms which appear to be physically relevant, but
  previously not discovered in the literature.}

{In this paper, we do not include the external magnetic
  field in the Schr\"odinger equation. Although the second order
  correction in the Bloch dynamics from the magnetic field results in
  many interesting applications, rigorous justification of such terms
  remains a challenge in mathematics, which we leave for future
  studies.}

The rest of the paper is organized in the following way. We present a brief review of the theory of Bloch decomposition and introduce the framework of the perturbation method to the Bloch wave function in Section 2. In Section 3, we carry out a systematic two-scaled WKB analysis  to the Schr\"odinger equation with a lattice potential and a slow-varying external potential, where the phase equation up to second order corrections has been derived and the validity of the WKB ansatz has been justified. At last, we show in Section 4  the  derivation of the characteristic equations of the phase equation with second order corrections, and under certain physical assumptions, the characteristic equation reduces to  the bi-characteristic equations with corrected Berry curvature.

Throughout this paper, we assume the following the convention in notation. If an $\epsilon$ dependent  function $f^\eps$ admits an asymptotic expansion, we denote the $n-$th order term by $f^n$, and the sum of the first $n+1$ terms by $f_{(n)}$, namely,
\[
f^\eps = f^0 + \eps f^1+\cdots +\eps^n f^n + \mathcal{O}(\e^{n+1})= f_{(n)} + \mathcal{O}(\e^{n+1}).
\]
Also we use notations as $A^{(n)}$ to stress that it is a $n-$th order
tensor.

\section{Preliminaries and the static perturbation} 


\subsection{Bloch decomposition}

Recall the Schr\"odinger equation with a periodic lattice potential
and a slow-varying scalar potential
\begin{equation}\label{SE}
i\eps \frac{ \partial}{\partial t}\psi(t,x)=H \psi(t,x)= \left(-\frac{\eps^2}{2}\Delta_x+V\left( \frac{x}{\eps}\right)+U(x)\right)\psi(t,x).
\end{equation}
In the absence of the external potential $U$, the Hamiltonian, after a
change of variable $z = x / \eps$, is given by
\[
H_{\text{per}} = -\frac{1}{2}\Delta_z + V(z).
\]
It is translational invariant with respect to the lattice
$\mathbb{L}$. 
As a result, the spectrum of the Hamiltonian can be
understood by the Bloch-Floquet theory, see e.g., \cite{RS}. In
particular, we have the periodic Bloch wave functions $\Psi^0_n(z,p)$,
given as the eigenfunctions of \footnote{The superscript $0$ stands
	for ``unperturbed'', as will be clear in the sequel.}
\begin{equation}\label{eigen01}
H^0(p) \Psi_n^0(z,p) :=\left( \frac{1}{2} \left(-i\nabla_z +p \right)^2 +V(z) \right)  \Psi_n^0(z,p)=E^0_n(p)\Psi^0_n(z,p) 
\end{equation}
on $\Gamma$ with periodic boundary conditions. Here $\Gamma$ is the
\emph{unit cell} of lattice $\mathbb{L}$ and $p \in \Gamma^{\ast}$ is
the crystal momentum, where $\Gamma^{\ast}$ denotes the \emph{first
	Brillouin zone} (unit cell of the reciprocal lattice). 
For each fixed $p\in \Gamma^*$, the Bloch Hamiltonian $H^0(p)$ is a
self-adjoint operator with compact resolvent, the spectrum of which is
given by
\begin{equation*} 
\sigma \left(H^0(p) \right) = \bigl\{E^0_n(p) \mid
n \in \Z_+ \bigr\} \subset \R, 
\end{equation*}
where the eigenvalues $E^0_n(p)$ (counting multiplicity) are
increasingly ordered $E^0_1(p) \le \cdots \le E^0_n(p) \le E^0_{n+1}(p) \le \cdots$.
It is shown by Nenciu \cite{Neff} that for any $n\in \Z_+$ there
exists a closed set $C_n \subset \Gamma^*$ of measure zero
such that
\[
E_{n-1}^0(p)<E_n^0(p)<E_{n+1}^0(p),\quad p\in \Gamma^* \backslash C_n, 
\]
and moreover $E_n^0(p)$ and $\Psi^0_n(\cdot,p)$ are analytic on $p\in
\Gamma^* \backslash C_n$.  In this paper, we stick to the adiabatic
regime and assume the energy band $n$ of interest is an isolated Bloch
band (\textit{i.e.}, $C_n = \emptyset$).  As a consequence, $E_n^0(p)$
and $\Psi_n^0(\cdot, p)$ are analytic with respect to $p$ in
$\Gamma^{\ast}$. Moreover, as we will focus on the particular single
band, we will suppress the subscript $n$ unless otherwise indicated.

Given $f \in L^2(\R^d)$, we recall the \emph{Bloch transform}, which is an isometry from $L^2(\R^d)$ to $L^2(\Gamma \times \Gamma^{\ast})$
\begin{equation}\label{psiphi}
\wt{f}(z, p) = \frac{\lvert\Gamma\rvert}{(2\pi)^d} \sum_{X \in \mathbb{L}} f(z+ X) e^{-ip\cdot(z+ X)},
\end{equation} 
where $\lvert\Gamma\rvert$ denotes the volume of the unit cell of the lattice $\mathbb{L}$. 
The inverse transform is given by 
\begin{equation}\label{phipsi}
f(z) = \int_{\Gamma^{\ast}} e^{i p \cdot z} \wt{f}(z, p) \,\mathrm{d} p. 
\end{equation}

\subsection{Perturbed Bloch wave functions}\label{PBloch}

With the external potential $U$, our analysis needs perturbation of
the Bloch waves. For this, let us recall the perturbation theory in
the context of Bloch wave functions.  We assume a family of
Hamiltonian $H^{\eps}(x, p)$ on $L_z^2(\Gamma)$ parametrized by $x$ and
$p$ admits the \emph{static} asymptotic expansion
\[
H^\eps(x, p) =H^0(p) +\e H^1(x,p)+\e^2H^2(x, p)+\mathcal{O}(\e^3). 
\] 
This expansion is called static because it does not capture the
dynamical information in the time propagation.  We assume the leading
order term to be just given by the Bloch Hamiltonian $H^0(p) =
\tfrac{1}{2} (-i \nabla_z + p)^2 + V(z)$, which is independent of $x$.
For the eigenvalue problem
\[
H^\e\Psi^\e=E_s^\e \Psi^\e, 
\]
we assume the asymptotic expansions
\begin{align*}
& E_s^\e(p,x)=E^0(p)+\eps E_s^1(p,x)+\e^2 E_s^2(p,x)+\mathcal{O}(\eps^3), \\
& 
\Psi^\e(z,p,x)=\Psi^0(z,p)+\eps \Psi^1(z,p,x)+\e^2 \Psi^2(z,p,x)+\mathcal{O}(\eps^3).
\end{align*}
Here, we call $E_s^\eps$ the static Bloch energy, in contrast to the dynamic Bloch energy which we shall define later. But, we suppress the subcript in $E^0$ since $E^0$ is well-defined by the unperturbed Hamiltonian.
The leading (zeroth) order terms in $\eps$ yields
\begin{equation}
H^0(p) \Psi^0(\cdot, p) = E^0(p) \Psi^0(\cdot, p),
\end{equation}
which is just the unperturbed eigenvalue problem. In particular, $E_s^0$
and $\Psi^0$ are independent of the parameter $x$. On the other hand,
the higher order terms of the Hamiltonian will depend on $x$
explicitly, since these terms (to be specified later) will be used to
capture the inhomogeneous influence by the slow-varying scalar
potential $U(x)$. Thus, the higher order expansions of $E^{\eps}$ and
$\Psi^{\eps}$ depend on both $p$ and $x$.

Collecting the terms in the next order, one gets
\begin{equation}\label{eigen10}
H^0 \Psi^1+H^1 \Psi^0=E^0 \Psi^1 + E_s^1 \Psi^0.
\end{equation}
Taking the inner product with $\Psi^0$ and using the leading order
equation, one gets
\begin{equation}\label{for:E1}
E_s^1(p, x) =\left\langle \Psi^0(\cdot, p), H^1(p, x) \Psi^0(\cdot, p) \right\rangle.
\end{equation}
Rewrite the first order equation as
\[
(H^1-E_s^1) \Psi^0=(E^0-H^0)\Psi^1.
\]
Given $E_s^1$ and $H^1$, $\Psi^1$ is then determined up to an arbitrary
constant multiple of $\Psi^0$. To fix the arbitrariness, we will take
$\langle \Psi^1, \Psi^0 \rangle = 0$, which turns out to simplify some
of the calculations in our analysis.

The equation of the second order terms reads 
\begin{equation}\label{eigen20} 
H^0 \Psi^2 +H^1 \Psi^1 + H^2
\Psi^0 = E^0 \Psi^2 + E_s^1 \Psi^1 + E_s^2 \Psi^0.
\end{equation}
To get $E_s^2$, it suffices to take the inner product with $\Psi^0$
\begin{equation}\label{E2}
\begin{aligned}
E_s^2(p, x) & =\left\langle \Psi^0(\cdot, p), H^1(p, x)
\Psi^1(\cdot, p, x) \right\rangle - E_s^1(p, x)\left\langle
\Psi^0(\cdot, p), \Psi^1(\cdot,
p, x) \right\rangle \\
& \qquad \qquad + \left\langle \Psi^0(\cdot, p), H^2(p, x)
\Psi^0(\cdot, p) \right\rangle \\
& = \left\langle \Psi^0(\cdot, p), H^1(p, x) \Psi^1(\cdot, p, x)
\right\rangle + \left\langle \Psi^0(\cdot, p), H^2(p, x)
\Psi^0(\cdot, p) \right\rangle, 
\end{aligned}
\end{equation}
where the second equality follows from
$\langle \Psi^1, \Psi^0 \rangle = 0$.  Moreover, $\Psi^2$ can be
solved from \eqref{eigen20} given $E_s^2$. This procedure can be
continued to higher orders, which we will omit as only the corrections
up to the second order will be considered in this paper.

To derive explicit formulas for the static expansion, we consider the Hamiltonian in \eqref{SE}. To treat the slow-varying potential $U(x)$, we assume we consider $x$ in the neighborhood of a point $x_c$, with $x=x_c + \eps z$ satisfied. Then, we do Taylor expansion of $U(x)$ around $x_c$ and get
\[
U(x)=U(x_c+\eps z)=U(x_c)+ \eps z \cdot \nabla U(x_c) +\eps^2 \frac 1 2 z \cdot \nabla^2 U(x_c) z +\mathcal O(\eps^3). 
\]
Clearly, $\eps z \cdot \nabla U(x_c)$ corresponds the the first order correction to the unperturbed Hamiltonian, and $\eps^2 \frac 1 2 z \cdot \nabla^2 U(x_c) z$  accounts for the second order correction. We remark that, the perturbation in Hamiltonian in the form of $\eps z \cdot \nabla U(x_c)$ accounts for for many physics phenomenon, such as the Wannier-Stark ladders, see e.g., \cite{Neff}. 

To see how a perturbation like $z \cdot \nabla U(x_c)$ can be related to  crystal momentum $p$, we note that using  Bloch transformation, we have
\begin{align*}
z f(z) = \int_{\Gamma^{\ast}} z e^{ip\cdot z} \wt{f}(z, p) \,
\mathrm{d} p = \int_{\Gamma^{\ast}} -i \nabla_p e^{ip\cdot z}
\wt{f}(z, p) \, \mathrm{d} p = \int_{\Gamma^{\ast}} i e^{ip\cdot z}
\nabla_p \wt{f}(z, p) \, \mathrm{d} p.
\end{align*}
Hence, we arrived at  the first order correction to the Hamiltonian 
\[
H^1(p,x_c)= i  \nabla U(x_c) \cdot \nabla_p, 
\]
and according to \eqref{for:E1}, the first order correction to the static Bloch energy is
\begin{equation}\label{eq:E10}
E_s^1 (p,x_c)= \nabla U(x_c) \cdot i\langle \Psi^0, \nabla_p \Psi^0\rangle. 
\end{equation}
Similarly, by \eqref{E2}, we obtain the second order correction  to the static Bloch energy
\begin{equation}\label{eq:E20}
E_s^2(p,x_c)=\nabla U(x_c) \cdot i\langle \Psi^0, \nabla_p \Psi^1\rangle+\frac 1 2  \nabla^2 U(x_c) : \langle \Psi^0, \nabla_p^2 \Psi^0 \rangle.
\end{equation}
 {
We remark that the term $i\langle \Psi^0, \nabla_p \Psi^1\rangle$ in \eqref{eq:E20} is related to the first order correction to the Berry curvature, which we define in the next section. For convenience, we also define the real part of $E_s^2$
\begin{equation}\label{eq:E30}
	\tilde{E}_s^2=\frac 12\nabla U(x_c) \cdot \left( i\langle \Psi^0, \nabla_p \Psi^1\rangle+c.c.\right)+\frac 1 4  \nabla^2 U(x_c) : \left(\langle \Psi^0, \nabla_p^2 \Psi^0 \rangle+c.c.\right),
\end{equation}
which naturally shows up in the WKB analysis in Section \ref{secWKB}.

	\subsection{Deriving the first order correction of the Berry connection}
	As we have mentioned in Section \eqref{PBloch}, some important terms in the static Bloch energy are known as the Berry connections. Here, we give the definition of the Berry connection and we shall reformulate the first order Berry connection to an expression in terms of the leading order Bloch wave functions. In Section \eqref{secWKB}, we shall see that those terms are also involved in the second order Bloch dynamics. 

	Define the leading order Berry connection as
	\begin{equation}
		{\mathcal{A}}^0(p)=i\langle \Psi^0, \nabla_p \Psi^0 \rangle,
	\end{equation}
	and the first order correction of the  Berry connection as
	\begin{equation}\label{con13}
		{\mathcal{A}}^1(t,p,x)=\frac 1 2 \left( i \langle \Psi^0, \nabla_p \Psi^1 \rangle + c.c. \right).
	\end{equation} 
%
%
	We also introduce the velocity operator, which is also used by the physicists (see \cite{GYN}), 
	\begin{equation}\label{con3}
	\hat{V} = -i[z,H^0]=-i(zH^0-H^0z).
	\end{equation}
%
	It is easy to show that the operator $\hat{V}$ in equation \eqref{con3} is equivalent to $\nabla_pH^0$.  In fact,
	\begin{align*}
	\hat{V}f(z) &= -i(zH^0-H^0z)f(z), \\
	    &= -i(z(\frac{1}{2}(-i\nabla_z+p)^2+V(z))-(\frac{1}{2}(-i\nabla_z+p)^2+V(z))z)f(z),
	\end{align*}
	and after some simplification, we get
	

	\begin{displaymath}
	\hat{V}f(z) = (-i\nabla_z+p)f(z).
	\end{displaymath}
	On the other hand, from \eqref{eigen01} we have
	\begin{displaymath}
	\nabla_pH^0(p)=-i\nabla_z+p,
	\end{displaymath}
	and thus we conclude $\hat{V}=\nabla_pH^0$ in the operator sense.

	Next, by differentiating the equation \eqref{eigen01} with respect to $p$, we have
	\begin{equation}\label{con6}
	\nabla_pH^0\Psi_m^0 + H^0\nabla_p\Psi_m^0=\nabla_pE^0_m\Psi_m^0 + E^0_m\nabla_p\Psi_m^0.
	\end{equation}
	We replace $\nabla_pH^0$ with $\hat{V}$, take the inner product with $\Psi_n^0$ ($n\neq m$), and obtain
	\begin{equation}
	<\Psi_n^0, \hat{V}\Psi_m^0>+E_n^0<\Psi_n^0, \nabla_p\Psi_m^0>=E_m^0<\Psi_n^0, \nabla_p\Psi_m^0>.
	\end{equation}
	Now we have,
	\begin{equation}\label{con9}
	<\Psi_n^0, \nabla_p\Psi_m^0>=\frac{<\Psi_n^0, \hat{V}\Psi_m^0>}{E_m^0-E_n^0},
	\end{equation}
	or we can write it in the index form
	\begin{equation}
	<\Psi_n^0, \partial_{p_j}\Psi_m^0>=\frac{<\Psi_n^0, \hat{V}_j\Psi_m^0>}{E_m^0-E_n^0},
	\end{equation}
	where $\hat{V}_j$ is the j-th component of $\hat{V}$.

	For the rest of this section,  we assume that the external potential $U$ is linear, and thus we have $-\nabla U= \vec{E}$, where $\vec{E}$ is the constant electric field.  Please note that, the assumption is taken only to reach the agreement with the results in the physics literature. Then, the first order perturbation of the Hamiltonian is simplified to
	\begin{equation}
	H^1 =-i\vec{E}\cdot\nabla_p.
	\end{equation}


	Without loss of generality, we focus on the first order correction of the Berry connection associated with the first Bloch band. 
	Collecting the terms in the first order of the Bloch equation, 
	\begin{equation}\label{con14}
	H^0 \Psi_1^1+H^1\Psi^1_0=E_0 \Psi_1^1+E^1\Psi^1_0.
	\end{equation}
	Due to the completeness of the leading order Bloch functions, and the fact that  $\langle \Psi^0_1, \Psi_1^1\rangle =0$, we write 
	\begin{equation} \label{eq:phi11exp}
	\Psi_1^1=\sum_{k=2}^{\infty}C_k\Psi_k^0,
	\end{equation}
	and the equation \eqref{con14} becomes
	\begin{equation}
	\sum_{k=2}^{\infty}C_k E_k^0 \Psi_k^0-i\vec{E}\cdot\nabla_p\Psi^0_1=\sum_{k=2}^{\infty}C_kE_1^0\Psi_k^0+E^1_1\Psi_0^0.
	\end{equation}
	Take the inner product with $\Psi_k^0$ ($k\neq 1$), we obtain 
	\begin{equation}
	C_k=\frac{i\vec{E}\cdot<\Psi_k^0,\nabla_p\Psi_1^0>}{E_k^0-E_1^0}.
	\end{equation}
	
	
	With equation \eqref{con9}, we can further rewrite 
	\begin{equation} \label{eq:coef}
	C_k=\frac{i\vec{E}\cdot<\Psi_k^0,\hat{V}\Psi_1^0>}{(E_k^0-E_1^0)^2}=\sum_{j}\frac{i(\vec{E})_j<\Psi_k^0,\hat{V}_j\Psi_0^0>}{(E_k^0-E_1^0)^2}.
	\end{equation}
 Substitute \eqref{eq:phi11exp} and \eqref{eq:coef} into equation \eqref{con13}, we obtain the first order correction of the Berry connection as an expression of the leading order Bloch functions and Bloch energies, 
	\begin{equation}
	\mathcal{A}^1(t,p,x)={\rm Re}\left(-\sum_{k=2}^{\infty}\frac{\vec{E}\cdot<\Psi_k^0,\hat{V}\Psi_1^0><\Psi_1^0,\hat{V}\Psi_k^0>}{(E_k^0-E_1^0)^3}\right)
	\end{equation}
The equivalent index form is given by
	\begin{displaymath}
	\mathcal{A}^1_i(t,p,x)=\sum_{j}{\rm Re}\left(\sum_{k=2}^{\infty}\frac{(\vec{E})_j<\Psi_k^0,\hat{V}_j\Psi_1^0><\Psi_1^0,\hat{V}_i\Psi_k^0>}{(E_1^0-E_k^0)^3}\right)
	\end{displaymath}
	In fact, this result is standard from the perspective of the
    perturbation theory, which exactly agrees with equation (1) of
    \cite{GYN} in the absence of the external magnetic field. Whereas,
    investigating the role of the Berry connection with its first
    order correction in dynamics is more challenging.

    The  work \cite{GYN} proposed a derivation of the effective
    Lagrangian for the wave packet dynamics when the external fields
    are constant, and its Euler-Lagrangian equations with a positional
    shift lead to effective semi-classical dynamics with the second
    order corrections. In the next section, we aim to apply a
    two-scale WKB asymptotic analysis to obtain the Bloch dynamics
    with the second order corrections. With the same positional shift,
    we shall see our results essentially agree with the results in
    \cite{GYN} and thus provide a mathematical justification.

	
	
}

\section{WKB asymptotic analysis}\label{secWKB}
\subsection{The ansatz and the zeroth order equation}

In this section, we carry out a two-scaled mono-kinetic WKB analysis
to the Schr\"odinger equation.  The starting point is the following
ansatz, which is a natural extension of the ones applied in
\cite{ELY}, to the Schr\"odinger equation \eqref{SE}.
\begin{equation} \label{WKB}
\psi_{\rm w}  = A^{\eps}(t,x)  \chi^{\eps}\left(\frac{x}{\e},t,\nabla_x S^ \eps,x \right) \exp\left(\frac{i}{\eps}S^{\eps}(t,x)\right), 
\end{equation}
where $\chi^\e(z,t, p,x)$ is the modified Bloch waves with the
asymptotic expansion
\[
\chi^\e(z,t,p,x)=\chi^0(z,p)+ \eps \chi^1(z,t,p,x)+ \eps^2 \chi^2(z,t,p,x)+\mathcal O(\eps^3).
\]
We will take
\[
\chi^0(z,p)=\Psi^0 (z,p),
\]
which does not depend on $t$ or $x$, and we expect $\chi^k(z,t,p,x)$
contains $\Psi^k(z,p,x)$ and necessary modification terms to be
specified. We emphasize that, even though $\Psi^k(z,p,x)$ is
time-independent, $\chi^k (z,t,p,x)$ might be time-dependent due to
the modification terms. We also assume the asymptotic expansions for the phase and amplitude in the ansatz \eqref{WKB}
\begin{align*}
  & S^{\eps}(t,x)=S^0(t,x)+\eps S^1(t,x) + \eps^2 S^2(t,x) +\mathcal{O}(\eps^3), 
\intertext{and}
  &  A^{\eps}(t,x)=A^0(t,x)+\eps A^1(t,x) + \eps^2 A^2(t,x) + \mathcal{O}(\eps^3).
\end{align*}
We emphasize that here the phase function series $\{S^k\}$ and the
amplitude function series $\{A^k\}$ are \emph{real-valued}, while the
functions $\chi^{\eps}$ are complex-valued, all yet to be
determined. Note that, since we will focus on a particular band
throughout the analysis, we have suppressed the energy band subscript
$k$ to simplify the notation.  The validity and the accuracy of this
ansatz will be discussed later.
 
Note that, we can rewrite the ansatz in the following way
\[
\psi_{\rm w}= a^\e\left(t,\frac{x}{\e},x\right) \exp(i S^0(t,x)/\e),
\]
where
\[
a^\e\left(t,\frac{x}{\e},x\right)= A^{\eps}(t,x)\chi^{\eps}\left(\frac{x}{\e},t,\nabla_x S^ \eps,x \right) \exp\left(i(S^1+\e S^2+\cdots) \right), 
\]
is the total amplitude, which is clearly complex-valued.  We remark
that, the ansatz \eqref{WKB} is special in the sense that, the total
amplitude $a^\e$ depend on the $\chi^n$ in a restrictive way. To be
more specific, since $A^k$ are all real-valued, it contains no phase
information at all while all the phase information is contained in the
term $\exp\bigl(i(S^1+\e S^2+\cdots) \bigr)$.  This assumption is
essential to guarantee that the canonical variables have a unique
trajectory. This is also  why this ansatz is called mono-kinetic.

We also emphasize that, as our purpose is to derive Bloch dynamics and its corrections, we do not aim to find a general approximate solution to the Schr\"odinger equation up to all time, but rather a specific solution which describes the behavior of a wave packet propagating under \eqref{SE}. In particular, \eqref{WKB} poses restrictions on the initial condition, which we will discuss further below. 

\smallskip 

A straightforward calculation yields
\begin{align*} 
 e^{-iS^\eps/\eps}\partial_t \psi_{\rm w} &= \partial_t A^{\eps} \chi^{\eps} + A^\eps \partial_t \chi^\eps +\frac{i}{\eps}A^{\eps}\chi^{\eps} \partial_t S^\eps+ A^{\eps}\nabla_p \chi^{\eps}\cdot \nabla_x \partial_t S^{\eps}, \\
 e^{-iS^\eps/\eps} \nabla_x \psi_{\rm w} &= \nabla_x A^{\eps} \chi^{\eps}+ \frac{i}{\eps} A^{\eps}\chi^{\eps}\nabla_x S^\eps + A^{\eps} \nabla_x \chi^{\eps} + A^{\eps} \nabla_x^2 S^\eps \nabla_p \chi^{\eps}+ \frac{1}{\eps}A^{\eps}\nabla_z \chi^{\eps},  
\intertext{and}
\displaystyle e^{-iS^\eps/\eps}\Delta_x \psi_{\rm w} &=  \Delta_x A^{\eps} \chi^{\eps} + 2 \frac{i}{\eps} \nabla_x A^{\eps} \cdot \nabla_x S^\eps \chi^{\eps} + 2 \nabla_x A^{\eps} \cdot \nabla_x \chi^{\eps} + 2 \nabla_x A^{\eps} \cdot \nabla^2_x S^\eps \nabla_p \chi^{\eps}  \\
&\quad + 2 \frac{1}{\eps} \nabla_x A^{\eps} \cdot \nabla_z \chi^{\eps} +\frac{i}{\eps} A^{\eps}  \chi^{\eps} \Delta_x S^\eps  - \frac{1}{\eps^2} |\nabla_x S^\eps|^2 A^{\eps}\chi^{\eps}+ 2 \frac{i}{\eps}A^{\eps} \nabla_x S^\eps \cdot \nabla_x \chi^{\eps} \\
&\quad + 2\frac{i}{\eps}A^{\eps}\nabla_x S^\eps \cdot \nabla^2_x S^\eps \nabla_p \chi^{\eps} +2 \frac{i}{\eps^2} A^{\eps} \nabla_x S^\eps \cdot \nabla_z \chi^{\eps}+A^{\eps}\Delta_x \chi^{\eps}\\
&\quad +2 A^{\eps} \nabla_x^2 S^\eps \nabla_p \cdot \nabla_x \chi^{\eps} + 2 \frac{1}{\eps} \nabla_x \cdot \nabla_z \chi^{\eps} A^{\eps}+\nabla_p \chi^{\eps} \cdot \left( \nabla_x \cdot \nabla^2_x S^\eps \right) A^{\eps}   \\
&\quad +A^{\eps} \left( \nabla_x^2 S^\eps \nabla_p \right)^2 \chi^{\eps}+ 2\frac{1}{\eps} A^{\eps} \nabla_x^2 S^\eps \nabla_p\cdot \nabla_z \chi^{\eps} +\frac{1}{\eps^2}A^{\eps}\Delta_z \chi^{\eps}.
\end{align*}
After sorting the terms in order of $\eps$, we obtain 
\begin{align}
  i\eps e^{-iS^\eps/\eps}\partial_t \psi_{\rm w} & = T_0+\eps T_1, \label{eq:tdpsi} \\
e^{-iS^\eps/\eps}\left( -\frac{\eps^2}{2}\Delta_x +V\left(\frac{x}{\eps} \right)+U(x) \right) \psi_{\rm w} & =F_0+\eps F_1 + \eps^2 F_2, \label{eq:sdpsi}
\end{align}
where we have introduced the short-hand notations 
\begin{align*}
  & T_0= - A^{\eps}\chi^{\eps} \partial_t S^\eps, \quad T_1=i \partial_t A^{\eps} \chi^{\eps}+ i A^\eps \partial_t \chi^\eps+i A^{\eps}\nabla_p \chi^{\eps}\cdot \nabla_x \partial_t S^{\eps}, \\
  & F_0 = \frac{1}{2} |\nabla_x S^\eps|^2 A^{\eps}\chi^{\eps}- i A^{\eps} \nabla_x S^\eps \cdot \nabla_z \chi^{\eps}-\frac{1}{2}A^{\eps}\Delta_z \chi^{\eps}+V\left(\frac{x}{\eps} \right)A^{\eps}\chi^{\eps}+U(x)A^{\eps}\chi^{\eps} \\
  & \quad =\left( H^0(p) \left( \frac{x}{\eps},\nabla_x S^\eps
    \right)+U(x)\right)\chi^\eps A^ \eps, \\
  & F_1  = -\frac{i}{2} \left(\Delta_x S^\eps +2 (\nabla_x S^\eps -i\nabla_z)\cdot \nabla_x^2 S^\eps \nabla_p  \right)\chi^\eps A^ \eps-i(\nabla_x S^\eps -i\nabla_z)\cdot \nabla_x \chi^ \eps A^ \eps \\
  & \qquad -i \nabla_x A^0 \cdot (\nabla_x S^\eps -i\nabla_z)\chi^\eps. \\
  & F_2 =-\frac{1}{2} \Delta_x A^{\eps} \chi^{\eps} - \nabla_x A^{\eps} \cdot \nabla_x \chi^{\eps}- \nabla_x A^{\eps} \cdot \nabla^2_x S^{\eps} \nabla_p \chi^{\eps}-\frac{1}{2}A^{\eps}\Delta_x \chi^{\eps}- A^{\eps} \nabla_x^2 S^\eps \nabla_p \cdot \nabla_x \chi^{\eps}  \\
  & \qquad -\frac{1}{2}\nabla_p \chi^{\eps} \cdot \left( \nabla_x \cdot
    \nabla^2_x S^\eps \right) A^{\eps} -\frac{1}{2} A^{\eps} \left(
    \nabla_x^2 S^\eps \nabla_p \right)^2 \chi^{\eps}.
\end{align*} 
Combining \eqref{eq:tdpsi} and \eqref{eq:sdpsi} and matching by order
of $\eps$, to the leading order, we get the following Hamilton-Jacobi
equation for the leading term of the phase function
\begin{equation}\label{eqS0}
-\partial_t S^0 = E^0(\nabla_x S^0)+ U(x),
\end{equation}
where we have used the identity \eqref{eigen01}. Recall that, for an isolated Bloch band, $E^0(p)$ is analytic for all $p\in \Gamma^*$. Thus, the equation of $S^0$ can be solved by the method of characteristics, where the characteristic flow is determined by the following Hamiltonian equations: 
\[
 \dot Q= P,    \quad
 \dot P = -\nabla_P E^0(Q), 
\]
with initial conditions 
\[
Q(0)=x,\quad P(0)=\nabla_x S^0(0,x).
\]

We remark that, the characteristic lines obtained by the above
Hamiltonian flow are interpreted as the rays of geometric
optics. Given the initial phase $S^0(0, x)$, the Hamiltonian system
locally defines a flow map. Caustics may appear at some finite time
when the characteristics initiated at difent locations
intersect. In the event of caustics, the Hamiltonian system no long
has classical solutions and the WKB solutions to the Schr\"odinger
equation \eqref{SE} breaks down as well. But since we aim to derive
corrections to the phase equation and to Bloch dynamics, it suffices
to consider the the WKB solutions before caustics formation.

It is natural to interpret $Q$ and $P$ as canonical variables of the
Hamilton-Jacobi equation. In Bloch dynamics, one considers a localized
wave packet, in which state the expectations of the position operator
and the momentum operator are called the classical variables. In the
semiclassical limit, to the leading order, the classical variables
agree with the canonical variables (see e.g., \cites{XCN, PST}).  As we
will see in later sections, within the Bloch dynamical equations with
corrections there is position shift introduced by Berry connections.

\subsection{The first order corrections}

\subsubsection{Useful identities for Bloch waves}
To study the first order correction, we shall first derive some useful identities
from the leading order Bloch eigenvalue problem:
\[
H^0(p)\chi^0(z,p) = E^0(p) \chi^0(z,p).
\]
difentiating this equation with respect to $p$, we get
\begin{equation} \label{eigen02}
(p-i\nabla_z)\chi^0 + H^0(p) \nabla_p \chi^0 = \nabla_p E^0 \chi^0 + E^0 \nabla_p \chi^0.
\end{equation}
The inner product with $\chi^0$ gives (since $(H_p^0 - E^0(p)) \chi^0 = 0$)
\begin{equation} \label{eigen03}
\left\langle \chi^0, (p-i\nabla_z)\chi^0 \right\rangle = \nabla_p E^0.
\end{equation}
Let us difentiate \eqref{eigen02} once again with respect to $p$, take a product with $\nabla_x^2 S^\eps$ and sum over indices (i.e., a Frobenius inner product for matrices), we arrive at 
\begin{multline} \label{eigen04}
  \Delta_x S^\e \chi^0 + 2 (-i \nabla_z +p) \cdot \nabla_x^2 S^\eps \nabla_p \chi^0 + H^0(p) \nabla_x^2 S^\eps \nabla_p \cdot \nabla_p \chi^0 = \\
  \nabla_x^2 S^\eps \nabla_p \cdot \nabla_p E^0 \chi^0 + 2 \nabla_x^2 S^\eps
  \nabla_p E^0 \cdot \nabla_p \chi^0 + E^0 \nabla_x^2 S^\eps \nabla_p
  \cdot \nabla_p \chi^0,
\end{multline}
which implies after taking inner product with $\chi^0$,
\begin{equation} \label{eigen05}
 \left\langle \chi^0,\left(\Delta_x S^\e + 2 (-i \nabla_z +p) \cdot \nabla_x^2 S^\eps \nabla_p\right) \chi^0 \right \rangle= \nabla_x^2 S^\eps \nabla_p \cdot \nabla_p E^0 + 2\nabla_x^2 S^\eps \nabla_p E^0 \cdot \left \langle \chi^0, \nabla_p \chi^0 \right \rangle. 
\end{equation}
We shall use these identities in our asymptotic derivation. 

\subsubsection{Derivation of the phase equation with first order correction}
Now we collect $\mathcal{O}(1)$ and $\mathcal{O}(\eps)$  terms from \eqref{eq:tdpsi} and \eqref{eq:sdpsi} and get 
\begin{equation}\label{order1}
  \begin{aligned} 
    &-A^0 \partial_t S^\eps \chi^0 + \eps i\partial_t A \chi^0 -\eps A^1 \partial_t S^\eps \chi^0 -\eps A^0 \partial_t S^\eps \chi^1 +i \eps A^0 \nabla_p \chi^0 \cdot \nabla_x \partial_t S^\eps  \\
    &\qquad \qquad =\left( H^0(p) \left( \frac{x}{\eps},\nabla_x S^\eps \right)+U(x)\right)\chi^0 A^0 \\
    &\qquad \qquad  + \eps\left( H^0(p) \left( \frac{x}{\eps},\nabla_x S^\eps \right)+U(x)\right)\chi^1 A^0 +\eps\left( H^0(p) \left( \frac{x}{\eps},\nabla_x S^\eps \right)+U(x)\right)\chi^0 A^1 \\
    &\qquad \qquad -\frac{i\eps}{2} 
      \left(\Delta_x S^\eps +2 (\nabla_x S^\eps -i\nabla_z)\cdot
        \nabla_x^2 S^\eps \nabla_p \right)\chi^0  A^ 0  \\
  &\qquad \qquad  -i\eps \nabla_x A^0 \cdot (\nabla_x S^\eps
      -i\nabla_z)\chi^0.
  \end{aligned} 
\end{equation}
We take inner product of both hand sides with $\chi^0$ and simplify
the result using identities \eqref{eigen03} and \eqref{eigen05}, we
obtain
\[
i \e \partial_t A^0 -A^0 \partial_t S^\e-i \e A^0 \langle \chi^0, \nabla_q \chi^0 \rangle \cdot \nabla_x U = (E^0+U(x))A^0 -i \e \nabla_x A^0 \cdot \nabla_p E^0- \frac {i\e}{2}\nabla_x \cdot \nabla_p E^0 A^0. 
\]
By separating the real and imaginary parts of the above equation, we
get
\begin{align} \label{eqS1}
& \partial_t S_{(1)} + U(x)+\left[ E^0+ i \e \langle \chi^0, \nabla_p \chi^0 \rangle \cdot \nabla_x U \right]\big\vert_{p=\nabla_x S_{(1)}}=0, \\
\label{eqA0}
& \partial_t A^0 =-\nabla_x A^0 \cdot \nabla_p E^0- \frac{1}{2}\nabla_x \cdot \nabla_p E^0 A^0
\end{align}
where we denote by $S_{(1)}=S^0+\e S^1$.
Let us consider the structure of these two equations. The leading
order equation \eqref{eqA0} for the amplitude
function is a transport equation. The phase equation \eqref{eqS1} with
the first order correction is still a Hamilton-Jacobi type equation. If we define
\[
E^1(p,x)=i \langle \chi^0, \nabla_p \chi^0 \rangle (p) \cdot \nabla_x
U(x),\] then the correction term in the phase equation is clearly
$E^1(\nabla_x S_{(1)},x)$. This correction term agrees with exactly
with the first order correction to static Bloch energy, and we will
identity this term with the first order correction to Bloch energy in
a \emph{dynamic} picture in Section \ref{E1dev}.  If we further
difentiate \eqref{eqS1} with respect to $x$, we get
\begin{equation*}
\partial_t \nabla_x S_{(1)}+ \nabla_x^2 S_{(1)} \cdot \nabla_p \left( E^0+ \e E^1 \right) + \nabla_x U + i \eps \nabla_x^2 U \cdot \langle \chi^0, \nabla_p \chi^0 \rangle= 0. 
\end{equation*}
This implies the bi-characteristics in canonical variables with the
notation $P=\nabla_x S_{(1)}$
\begin{align}  \label{eqX}
\dot Q &=  \nabla_p \left( E^0+ \e E^1 \right)|_{p=P,x=Q}\,  \\
 \label{eqP}
\dot P &= -  \nabla_x U (Q) - i \eps \nabla_x^2 U \cdot \langle \chi^0, \nabla_p \chi^0 \rangle|_{p=P,x=Q}.
\end{align} 
Here, $i\langle \chi^0, \nabla_p \chi^0 \rangle={\mathcal{A}}(p)$ is
known as the Berry connection or Berry vector potential. This quantity
is gauge-dependent, which means if one chooses difent phase factors
for $\chi^0$, the resulting ${\mathcal{A}}(p)$ are actually
difent. Namely, for an arbitrary smooth function $\zeta(p)$, if the
so-called gauge transformation
\[
\chi^0 \rightarrow \exp(i \zeta(p))\chi^0,
\]
is performed, the corresponding change happens for the Berry connection,
\[
{\mathcal{A}}(p) \rightarrow {\mathcal{A}}(p) - \nabla_p \zeta(p). 
\]
By a change of variables that takes into account the position shift
between classical variables and canonical variables by the Berry
connection,
\[
Q=\widetilde Q - \e {\mathcal{A}} (\widetilde P), \quad P=\widetilde P,
\]
we obtain, after dropping the tilde,
\begin{align*} 
\dot Q &= \nabla_P E^0(P)+\e \nabla_Q U \times \nabla_P\times {\mathcal{A}}(P),  \\
\dot P &= -\nabla_Q U(Q) .
\end{align*} 
Here, $\nabla_p \times {\mathcal{A}}(p)$ is the Berry curvature, which
is gauge independent.  This implies, in the corrected Bloch dynamics,
an anomalous velocity is introduced by the response of Bloch electrons
to the external electric field. In other words, the characteristic
speed has been modified due to the correction term $i \e \langle
\chi^0, \nabla_p \chi^0 \rangle \cdot \nabla_x U $ in the phase
equation. We remark that, so far our results agree with previous works
on first order corrections to the Bloch dynamics, see \cites{ELY, PST,
	XCN}. The focus of this paper is to extend this to second order
corrections to the Bloch dynamics.

\subsubsection{Derivation of $H^1$ and $E^1$}\label{E1dev}

In this part, we aim to derive the specific expressions of $H^1(p,x)$
and $E^1(p,x)$, keeping in mind that they should satisfy the equation
\eqref{SE} and the solution ansatz \eqref{WKB}.  Also, we will
establish the relation between $\chi^1$ with $\Psi^0$ and $\Psi^1$.


We observe that, since equation \eqref{eqS1} and equation \eqref{eqA0} have been derived, there is a difent perspective to view equation \eqref{order1}.
Substitute \eqref{eqS1} and \eqref{eqA0} into \eqref{order1}, we
obtain 
\begin{multline}\label{reorder1} 
{\mathcal{A}} \cdot \nabla_x U\chi^0 A^0 - i \nabla_p \chi^0 \cdot \nabla_x U A^0= \\
  + i \nabla_x A^0 \cdot (H^0(p)-E^0)\nabla_p \chi^0 + \frac{i}{2}
  (H^0(p)-E^0) \nabla_x^2 S^\eps \nabla_p \cdot \nabla_p \chi^0 A^0 +
  (H^0(p)-E^0)\chi^1 A^0. 
\end{multline}
The main idea here is to view this equation as the first order  perturbation equation for the \emph{dynamic} Bloch eigenvalue problem, as it has the same structure as  \eqref{eigen10}, in particular, all the terms with time derivatives in \eqref{order1} are canceled.
 
We decompose the perturbation $\chi^1$ into $\chi^1=w+v$, where
\begin{equation}\label{def_v}
  v= -\frac{i}{2} \nabla_x^2 S^0 \nabla_p \cdot \nabla_p \chi^0 - i \nabla_x \log A^0 \cdot \nabla_p \chi^0. 
\end{equation}
Unlike $v$, $w$ contains terms that cannot be written explicitly in $\chi^0$, which is given by 
\begin{equation}\label{eigen11}
{\mathcal{A}} \cdot \nabla_x U \chi^0  - i \nabla_p \chi^0 \cdot \nabla_x U = (H^0-E^0)w. 
\end{equation}
Note that this has the same structure of \eqref{eigen10} (recalled
here for convenience)
\begin{equation*}
 E^1 \Psi^0 -H^1 \Psi^0=(H^0-E^0)\Psi^1,
\end{equation*}
if we identify 
\[
\Psi^0=\chi^0,\quad \Psi^1=w,\quad E^1(p,x)={\mathcal{A}}(p) \cdot
\nabla_x U(x),\quad \text{and} \quad H^1(p, x) f= i \nabla_x U(x)
\cdot \nabla_p f.
\]
Note here $H^1$ is the first order correction to the unperturbed
Hamiltonian.  Clearly, the first order corrections to the Bloch energy
in the static expansion is of the same form as the first order
correction to the phase equation, and without confusion, hereafter, we
use $E^1$ to denote both meanings. However, in the dynamic picture,
the eigenfunction $\chi^1$ picks up an extra time-dependent part $v$
besides $\Psi^1$. This will have impact on the next
  order correction, as will be discussed below.

As we discussed in Section
\ref{PBloch}, it is natural to enforce the orthogonality condition
$\langle \chi^0, w \rangle=\langle \Psi^0,\Psi^1 \rangle=0$.  Moreover, this
condition will simplify the results in Section \ref{moreID}.  


We remark that $\chi^1$ depend on $x$ as a parameter, as clearly seen
from the definition of $v$ in \eqref{def_v} and the equation
\eqref{eigen11} for $w$. In addition, although $w$ does not depend on
$t$, $v$ depends on $t$ through $S^0$ and $\log A^0$, and we have
\[
\partial_t v(z,t,p,x) =-\frac{i}{2} \nabla_x^2 \partial_t S^0 \nabla_p \cdot \nabla_p \chi^0 - i \nabla_x \partial_t \log A^0 \cdot \nabla_p \chi^0,
\]
and hence $\chi^1$ depends on $t$.

\subsection{The second order corrections}
\subsubsection{More useful identities}\label{moreID}

For the second order correction, we need some more identities for the Bloch waves, following similar strategies as in Section 3.2.1 applied on \eqref{eigen11}. We conclude with the the following identities, we omit the details for the straightforward derivations.
\begin{equation}\label{eigen12}
\langle \chi^0, (-i \nabla_z +p) w \rangle + \langle \chi^0, H^1(p)\nabla_p \chi^0 \rangle = \nabla_p E^1+E^1 \langle \chi^0, \nabla_p \chi^0 \rangle,  
\end{equation}
\begin{multline} \label{eigen13}
\left\langle \chi^0,\left(\Delta_x S + 2 (-i \nabla_z +p) \cdot \nabla_x^2 S^\eps \nabla_p\right) w \right \rangle + \left\langle \chi^0, H^1 \nabla_x^2 S^\eps \nabla_p \cdot \nabla_p \chi^0 \right \rangle \\
 = 2 \nabla_x^2 S^\eps \nabla_p E^0 \cdot \langle \chi^0, \nabla_p w \rangle+ \nabla_x^2 S^\eps \nabla_p \cdot \nabla_p E^1   \\
+ 2 \nabla_p E^1 \cdot \langle \chi^0, \nabla_x^2 S^\eps \nabla_p \chi^0 \rangle +E^1 \langle \chi^0, \nabla_x^2 S^\eps \nabla_p\cdot \nabla_p \chi^0 \rangle,
\end{multline}
and
\begin{equation} \label{eigen14}
\langle \chi^0, (-i \nabla_z +p)\cdot \nabla_x w \rangle+ \langle \chi^0, \nabla_x \cdot (H^1 \nabla_p) \chi^0 \rangle = \nabla_x \cdot \nabla_p E^1 + \nabla_x E^1 \cdot \langle \chi^0, \nabla_p \chi^0 \rangle.
\end{equation}

Besides, from the definition of $v$, one gets
 \[
 \frac{i}{2} \nabla_x^2 S^\eps \nabla_p \cdot \nabla_p \chi^0= -v - i \nabla_x \log A^0 \cdot \nabla_p \chi^0.
 \] 
Then, by direct substitution and simplification, we obtain the following two identities,
\begin{align}\label{eigen15}
\frac{i}{2}A^0 E^1  \langle \chi^0, \nabla_x^2 S^\eps \nabla_p\cdot \nabla_p \chi^0 \rangle &= -A^0 E^1 \langle \chi^0, v \rangle -i E^1 \nabla_x A^0 \cdot  \langle \chi^0, \nabla_p \chi^0 \rangle.
\\
\label{eigen16}
-\frac{i}{2}A^0 \left\langle \chi^0, H^1 \nabla_x^2 S^\eps \nabla_p \cdot \nabla_p \chi^0 \right \rangle &= iA^0 \langle \chi^0,\nabla_p v\rangle \cdot \nabla_x U-\nabla_x A^0 \cdot \langle \chi^0, \nabla_p^2 \chi^0 \rangle \cdot \nabla_x U.
\end{align}

\subsubsection{Derivation of the corrected phase equation}\label{S2dev}

Now we collect the terms from \eqref{eq:tdpsi} and \eqref{eq:sdpsi} up
to $\mathcal{O}(\eps^2)$, this gives 
\begin{align*}
    &i\e \partial_t A^0 \chi^0 + i \e^2 \partial_t A^1 \chi^0 + i \e^2 \partial_t A^0 \chi^1 +i \eps^2 \partial_t v A^0 \\
    & \qquad - A^0 \partial_t S^\e \chi^0 -\e A^1 \partial_t S^\e \chi^0 - \e A^0 \partial_t S^\e \chi^1 -\e^2 A^0 \partial_t S^\e \chi^2  -\e^2 A^1 \partial_t S^\e \chi^1 - \e^2 A^2 \partial_t S^\e \chi^0  \\
    & \qquad + i \e A^0 \nabla_p \chi^0 \cdot \nabla_x  \partial_t S^\e +  i \e^2 A^1 \nabla_p \chi^0 \cdot \nabla_x  \partial_t S^\e+  i \e^2 A^0 \nabla_p \chi^1 \cdot \nabla_x  \partial_t S^\e  \\
    & = A^0 (H^0(p)+U)\chi^0 + \e A^1 (H^0(p)+U)\chi^0 + \e A^0 (H^0(p)+U)\chi^1   \\
    & \qquad + \e^2 A^2 (H^0(p)+U)\chi^0 + \e^2 A^1 (H^0(p)+U)\chi^1 + \e^2 A^0 (H^0(p)+U)\chi^2  \\
    &\qquad -i \e \nabla_x A^0 \cdot (\nabla_x S^\e -i \nabla_z) \chi^0 - i \e^2 \nabla_x A^1 \cdot (\nabla_x S^\e -i \nabla_z) \chi^0 -i \e^2  \nabla_x A^0 \cdot (\nabla_x S^\e -i \nabla_z) \chi^1  \\
    &\qquad -\frac{i\e}{2} A^0  \left( \Delta_x S^\e + 2 (-i \nabla_z +\nabla_x S^ \e)\cdot \nabla_x^2 S^\eps \nabla_p \right) \chi^0     \\
    &\qquad -\frac{i\e^2}{2} A^1  \left( \Delta_x S^\e + 2 (-i \nabla_z +\nabla_x S^ \e)\cdot \nabla_x^2 S^\eps \nabla_p \right) \chi^0   \\
    & \qquad -\frac{i\e^2}{2} A^0 \left( \Delta_x S^\e + 2 (-i \nabla_z +\nabla_x S^ \e)\cdot \nabla_x^2 S^\eps \nabla_p \right) \chi^1   \\
    & \qquad -i\e^2 (-i \nabla_z +\nabla_x S^ \e) \cdot \nabla_x \chi^1 A^0 -\frac{\eps^2}{2} \Delta_x A^0 \chi^0 -\eps^2 \nabla_x A^0 \cdot \nabla_x^2S^\e \nabla_p \chi^0  \\
    & \qquad -\frac{\e^2}{2}A^0 (\nabla_x \cdot \nabla_x^2 S^\e) \cdot \nabla_p \chi^0 -\frac{\eps^2}{2}A^0 (\nabla_x^2 S^\e \nabla_p)^2 \chi^0.	 
\end{align*}
Take the inner product with $\chi^0$ of the above equation, the real
part of the resulting identity reveals\footnote{The imaginary part
  leads to the equation that $A_1$ satisfies, which does not
  contribute to the second order corrections to Bloch dynamics, and
  will hence be neglected in this paper.}
\begin{equation}\label{eqS2}
\partial_t S_{(2)} +U(x) + \left[E^0 + \e E^1 + \e^2  E^2 \right]|_{p=\nabla_x S_{(2)}}  = 0,
\end{equation}
where $E^0$ and $E^1$ are of the form as before but are evaluated at
$\nabla_x S_{(2)}$ instead, and we denote by
$S_{(2)}=S^0+\e S^1+\e^2 S^2$. Here $E^2$ is the second order correction term to the phase equation.  Using the identities we derived before,
we find that $E^2$ takes the following form
\begin{equation}\label{eqE2}
\begin{aligned}
 E^2(p,x) &= \frac 12\left(i \langle \chi^0, \nabla_p w \rangle+c.c.\right) \cdot \nabla_x U +  \langle \nabla_p  \chi^0, \nabla_x^2 U \nabla_p \chi^0 \rangle - \frac{\Delta_x A^0}{2A^0} + \frac{1}{2} |\nabla_x^2 S^0 \nabla_p \chi^0|^2  \\
& \quad + \nabla_x \log A^0 \cdot {\rm Im} \left( \left\langle \chi^0, \left(p -i \nabla_z - \nabla_p E^0\right)v \right\rangle \right) - \frac{1}{2}\nabla_x^2 S^0 \nabla_p \cdot \nabla_p E^0 {\rm Im} \left( \langle \chi^0,v \rangle \right)  \\
& \quad -\nabla_x^2 S^0 \nabla_p E^0 \cdot {\rm Im} \left( \langle \chi^0, \nabla_p v \rangle \right) + \frac{1}{2} {\rm Im} \left( \left\langle  \chi^0, (\Delta_x S^0 +2 (p-i\nabla_z)\cdot \nabla_x^2 S^0 \nabla_p)v \right \rangle \right) \\
& \quad + {\rm Im} \left( \left\langle \chi^0, \left(p -i \nabla_z \right)\cdot \nabla_x v \right\rangle \right)+ {\rm Im} \left( \langle \chi^0,\partial_t v  \rangle \right). 
\end{aligned}
\end{equation}
Note that, here
\[
\nabla_x v = \nabla_x v(z,t,p,x)|_{z=\frac{x}{\e},\, p=\nabla_x S_{2}}.
\]
Since $\chi^0$ does not depend on $x$, we have
\[
\nabla_x v = -\frac{i}{2} \nabla_x^3 S^0  : \nabla_p^2  \chi^0 - i \nabla^2_x \log A^0 \cdot \nabla_p \chi^0.
\]
And by using equation \eqref{eqS0} and \eqref{eqA0}, we can directly compute that
\begin{align*}
\langle \chi^0, \partial_t v \rangle &= -\frac{i}{2} \nabla_x^2 \partial_t S^0: \left\langle \chi^0, \nabla_p^2 \chi^0 \right \rangle - i \nabla_x \partial_t \log A^0 \cdot \left\langle \chi^0, \nabla_p \chi^0 \right \rangle, \\
& =\frac{i}{2} \left( \nabla_x^3 S^0\cdot \nabla_p E^0+ (\nabla_x^2 S^0 \cdot \nabla_p)\otimes(\nabla^2_x S^0 \cdot \nabla)E^0 +\nabla_x^2 U \right) : \left\langle \chi^0, \nabla_p^2 \chi^0 \right \rangle \\
& + i \left( \nabla_x^2 \log A^0 \cdot \nabla_p E^0 + \nabla_x \log A^0 (\nabla_x^2 S^0 : \nabla^2_p E^0)\right) \cdot \left\langle \chi^0, \nabla_p \chi^0 \right \rangle \\
& + i \left( \frac{1}{2} \nabla_x^3 S^0 \cdot \nabla_p \cdot \nabla_p \cdot E^0+ \frac{1}{2} (\nabla_x^2 S^0 \cdot \nabla_p \cdot \nabla_p)(\nabla^2_x S^0 \cdot \nabla_p)E^0 \right) \cdot \left\langle \chi^0, \nabla_p \chi^0 \right \rangle.
\end{align*}
Note that, the last two lines are actually real-valued,  they do not enter the corrected phase equation, so we have
\[
{\rm Im}(\langle \chi^0, \partial_t v \rangle)= \frac{1}{2} \left( \nabla_x^3 S^0\cdot \nabla_p E^0+ (\nabla_x^2 S^0 \cdot \nabla_p)\otimes(\nabla^2_x S^0 \cdot \nabla)E^0 +\nabla_x^2 U \right) : {\rm Re} \bigl\langle \chi^0, \nabla_p^2 \chi^0 \bigr \rangle.
\]

Note, in the leading order phase equation, the unperturbed Bloch
energy $E^0$ naturally shows up, and the first order correction term
to the phase equation happens to be the first order correction of the
static Bloch energy $E^1$.  However, the second order correction to
the phase equation $ E^2$ does not necessarily agree with the second
order correction of the static Bloch energy $E_s^2$. Due to the WKB
ansatz \eqref{WKB} we have used, the phase equation and the amplitude
equation are no longer decoupled with second order corrections, and
thus some amplitude dependent terms, such as $\frac{\Delta_x
  A^0}{2A^0}$, enter the phase equation with second order
corrections. Besides, in the first order correction analysis, we have
learned that it is necessary to incorporate the term $v$ in $\chi^1$,
$v$ related terms also contribute to the second order correction to
the phase function.

 Recall that we have derived the second order correction of the Bloch energy \eqref{eq:E30} in the static expansion, repeated here for
convenience
\begin{equation}
\begin{aligned}
\tilde{E}_s^2  &= \frac 12\left(i \langle \chi^0, \nabla_p w \rangle+c.c.\right) \cdot \nabla_x U +  \frac 1 2 \langle \nabla_p  \chi^0, \nabla_x^2 U \nabla_p \chi^0 \rangle \\
& =\frac 12\left(i \langle \chi^0, \nabla_p w \rangle+c.c.\right) \cdot \nabla_x U - \frac{1}{4} \left( \langle \chi^0, \nabla_p^2 \chi^0 \rangle +\mbox{c.c.} \right): \nabla_x^2 U. 
\end{aligned}
\end{equation}
As we have mentioned, this is not necessarily the same as
${E}^2$ defined in \eqref{eqE2}, but ${E}^2$ contains all the terms in $\tilde{E}_s^2$ \eqref{eq:E30}.

\subsubsection{Well-posedness of initial data and accuracy of WKB approximation}
Before we turn to the characteristic equations from the phase equation
\eqref{eqS2}, let us summarize this section by discussing the
necessary assumptions such that the ansatz \eqref{WKB} is a valid
approximation and the corresponding accuracy before the formation of
caustics.

To this end, we reformulate \eqref{WKB} in the following form in order
to apply the convergence results in \cite{CMS} for two-scaled WKB
analysis.
\begin{align}\label{match}
\psi_{\rm w}&=\sum_{j=0}^{\infty}\e^j g_j (t, \frac{x}{\e},x)\exp\left(i \sum_{m=1}^{\infty} \e^{m-1} S^m(t,x) \right) \exp\left(\frac{i}{\eps}S^{0}(t,x)\right) \\
&:=\sum_{j=0}^{\infty} \e^j u_j (t, \frac{x}{\e},x)\exp\left(\frac{i}{\eps}S^{0}(t,x)\right), \nonumber
\end{align}
where
\[
u_j =g_j \exp\left(i \sum_{m=1}^{\infty} \e^{m-1} S^m(t,x) \right),
 \]
 and $g_j$ is determined by asymptotically matching the terms in \eqref{WKB}. For example,
\begin{equation}\label{g0}
g_0 = A^0(t,x) \chi^0(t,x/\e,\nabla_x S^0(t,x)),
\end{equation}
and
\begin{equation}\label{g1}
\begin{aligned}
  g_1&= A^0(t,x) \chi^1(t,x,x/\e,\nabla_x S^0(t,x)) + A^1(t,x) \chi^0(t,x/\e,\nabla_x S^0(t,x)) \\
  &\qquad +A^0(t,x) \nabla_x S^1(t,x)\cdot \nabla_p
  \chi^0(t,x/\e,\nabla_x S^0(t,x)).
\end{aligned}
\end{equation}

Taking $t = 0$, the ansatz above
imposes well-preparedness requirement on the initial condition to
guarantee the accuracy of this approximation. For simplicity, we
denote the truncated $J-$th order WKB approximation by
\[
\psi_{\rm w}^J:=\sum_{j=0}^{J-1} \e^j u_j \left(t, \frac{x}{\e},x\right)\exp\left(\frac{i}{\eps}S^{0}(t,x)\right).
\] 

To make the leading order WKB approximation valid, we make the
following assumptions, similar to those used in \cite{CMS}.
\begin{hyp}
  The initial wave function $\psi^\e_I$ is in the Schwartz space
  ${\mathcal{S}} (\R^d)$, and is of the WKB-type, \textit{i.e.}
  \[
  \psi_I^\e(x)=u_I\left(x,\frac{x}{\e}\right) e^{i\phi_I(x)/\e}+\eps
  \varphi_I^\e(x),
  \]
  with $\phi_I \in C^{\infty}(\R^d, \R)$,
  $u_I \in {\mathcal{S}} (\R^d\times \Gamma; \C )$. The function
  $\varphi_I^\e$ is to be specified later. The amplitude $u_I(x,y)$ is
  assumed to be concentrated on a single isolated Bloch band $E_k(p)$
  corresponding to a simple eigenvalue of $H^0$, i.e.
  \[
  u_I(x,y) = a_I(x) \chi_k(y, \nabla_x \phi_I(x)),
  \]
  where $a_I \in {\mathcal{S}}(\R^d;\C)$ is a given amplitude.
\end{hyp}
 One can check that, with the assumptions above, relation \eqref{g0} is always satisfied at $t=0$, and when $t=0$, $A^0$ and $S^0$ and $S^1$ are uniquely determined consequently. Moreover, from equation \eqref{reorder1}, $\chi^1$ is also determined at initial time. 
 
 Next, recall that we focus on wave functions which oscillate at the scale of $\mathcal{O}(\e)$, so we define the following spaces: for $s\in \N$, 
 \[
 ||f^\e||_{X^s_\e}=\sum_{|\alpha|+|\beta|\le s}||x^\alpha (\e \partial)^\beta f^\e||_{L^2},
 \]
 and define $X^s_{\e}$ as:
 \[
 X^s_\e=\left\{f^\e \in L^2 (\R^d); \sup_{0<\e\le1} ||f^\e||_{X^s_\e} \right\}.
 \]
 Now, we are ready to state the additional assumption for the next order WKB approximation.
\begin{hyp}
  (Well-prepared initial data.) The initial conditions $\psi_I^\e (x)$
   satisfy the assumptions above, and the
  leading part of the perturbation $\varphi_I^\e$ is of the particular
  form, so that there exist solutions to $A^1\in \R$ and $S^2\in\R$
  when $t=0$, whereas the residual part of $\varphi_I^\e$ is
  $\mathcal{O}(\e)$ in $X^s_\e$ for all $s\in \N$.
\end{hyp}
 
This assumption also implies that, initially each term in the
asymptotic expansion \eqref{WKB} is uniformly bounded in $X^s_\e$ for
all $s \in N$.
 
Note that, there are always initial conditions which satifiy all the
assumptions above, which are WKB-type initial condition concentrated
on one band with no tails. Then, by \cite{CMS}*{Theorem 4.5}, we obtain
the second order approximation to the exact solution.
\begin{theorem}
  Define $\Psi^\e(t)$ to be the solution to equation \eqref{SE} with
  initial conditions satisfying all the assumptions above. Assume
  there is no caustic formed before time $t_0$, the second order WKB
  approximation $\psi^2_W$ is valid up to any $t<t_0$, and for all
  $s\in \N$ there exists a constant $C$ such that
  \[
  \sup_{0<\tau<t} \lVert \Psi^\e(\tau)-\psi^2_w(\tau) \rVert_{X^s_\e}  \leq C\e^2.
  \] 
\end{theorem}

\section{Characteristic equations with second order corrections}
\label{sec:charac}

 
{ 
In this section, we derive the bi-characteristic equations for the
phase equation with the second order corrections. 
Recall that, the phase equation with the first order correction
\eqref{eqS1} is still a Hamilton-Jacobi equation, whose
characteristic equation in canonical variables $Q$ and $P$ are given
by
\begin{align*} 
\dot Q &=  \nabla_p \left( E^0+ \e E^1 \right)|_{p=P,x=Q},  \\
\dot P &= -  \nabla_x U (x) - i \eps \nabla_x^2 U \cdot \langle \chi^0, \nabla_p \chi^0 \rangle|_{p=P,x=Q}. 
\end{align*} 
However, for the phase equation with the second order corrections
\eqref{eqS2}, high order derivatives of $S^0$ as well as derivatives
of $\log A^0 $ are involved. This means, the corrected phase equation
does not admit a simple bi-characteristic structure. In other words,
such a characteristic system is not yet closed.

Whereas, one observes that the higher order derivatives of $S^0$ and
$\log A^0 $ are only contained in the terms of order
$\mathcal{O}(\e^2)$. Therefore, we can proceed in two ways:
\begin{enumerate}
\item One can solve the leading order phase  equation and
  amplitude equation in the pre-processing stage, then in the phase
  equation with second order corrections, all the higher order
  derivatives in the $\mathcal{O}(\e^2)$ terms are thus treated as
  known quantities. We take this treatment for the rest of the paper, 
and the physical meaning of the resulting characteristic equations are discussed in
  \S\ref{sec:disc}.
\item One can close the system of trajectories by introducing the
  characteristic equations of those higher order derivatives appeared
  in the $\mathcal{O}(\eps^2)$ correction. Since the results by 
 this approach are lengthy but straightforward, we shall omit the details in this work. 
\end{enumerate}
}

{ 
\subsection{The bi-characteristic equations in canonical variables}


Following \cites{GYN, GYN2}, we call
$ E^2$ the second order wave packet energy and define
\[
E_w(Q,P)= E^2(Q,P)-\tilde{E}_s^2(Q,P),
\]
where $E_w$ is interpreted as the extra wave packet energy besides the
static Bloch energy due to the specific profile of thel wave
function. In other words, $ E^2$ consists of two parts, the second
order correction to the static Bloch energy, and extra wave packet
energy which may be time-dependent due to its dependence on the phase
and the amplitude.  {We focus on the correction of the
  static Bloch energy part $\tilde{E}_s^2(Q,P)$ for the rest of this
  section. In particular, we prove that the bi-characteristic equation
  with $\tilde{E}_s^2(Q,P)$ is gauge invariant and analyze the
  physical interpretation of the bi-characteristic equations.}

{
Let us write 
 \[
 \tilde{E}_s^2 =  {\mathcal{A}}^1 \cdot \nabla_x U + {\mathcal{B}}:\nabla_x^2 U
 \]
 with the following notations}
 \[
 {\mathcal{B}}=-\frac{1}{4} \left( \langle \chi^0, \nabla_p^2 \chi^0 \rangle+ c.c. \right),
 \]
 and 
 \[{\mathcal{A}_{(1)}}={\mathcal{A}}^0+\e{\mathcal{A}}^1.\]
We recall that ${\mathcal{A}_{(1)}}$ is
referred to as the Berry connection with the first order correction which responds to the first derivative of the scalar potential. While
${\mathcal{B}}$ is a new quantity which is yet to be
explored. It responds to the second order derivative of the scalar potential and if only a uniform electric field is considered
as a perturbation to periodic Hamiltonians as in \cites{GYN, GYN2},
{$\tilde{E}^2_s$} loses the contribution from ${\mathcal{B}}$ because
$\nabla_x^2 U$ vanishes.

Then, we can write the bi-characteristic equations with the second order corrections as
\begin{align*}
\dot Q &= \nabla_P E^0 + \e \nabla_P ({\mathcal{A}}^0 \cdot \nabla_Q U) + \e^2 \nabla_P \tilde{E}_s^2(Q,P) + \e^2 \nabla_P E_w , \\
\dot P &= - \nabla_Q U -\e \nabla_Q^2 U \cdot {\mathcal{A}}^0 -\e^2 \nabla_Q \tilde{E}_s^2(Q,P) - \e^2 \nabla_Q E_w .
\end{align*}
We observe that, given that the auxiliary quantities involved in $E_w$ are known; the canonical coordinates still satisfy a Hamiltonian flow under a modified Hamiltonian
\[
\widetilde H_{(2)} =E^0(P)+ U(Q) +\e E^1(Q,P) + \e^2 \tilde{E}_s^2(Q,P)+ \e^2 E_w(Q,P).
\]
 We conclude this session by remarking that the bi-characteristic equations in canonical variables are valid for generic potential functions, whereas the correction terms might be gauge dependent. In the next, we proceed by considering special potential functions to investigate the bi-characteristic equations in physical variables.

\subsection{Physical interpretation of the bi-characteristic equations}
\label{sec:disc}

To convert the characteristic equations in canonical variables to
those in physical variable, one can apply a physics-motivated change
of variables to guarantee that the correction terms in new variables
are gauge invariant and thus physically meaningful. Such change of
variables are far from being fully understood in the second order
theory. When the external electrical field is uniform, the first order
correction of the Berry connection introduces an additional positional
shift in the change of variables (see \cite{GYN, GYN2}). We choose to
take such change of variables, and discuss its applications for two
scenarios: when the external electric field is uniform in space and
when the field is linearly varying in space.

  \subsubsection{Uniform electric field}

 To compare our results with those in recent work \cites{GYN, GYN2}, we follow their setup to  consider the special case taken in these papers where the electric
 field is uniform, or equivalently the scaler potential is linear. 
In this scenario, 
 $\nabla_x U$ reduces to a constant vector, and all the
 higher order derivatives of $U$ vanish. Moreover, $H^1$ and $E^1$ no
 longer depend on $x$, and as a result, ${\mathcal{A}}'={\mathcal{A}}^0
 + \e {\mathcal{A}}^1$ becomes $x-$independent.  
 
 Next, 
 we carry out the change of variables 
 that incorporates into the position shift due to the Berry connection, 
 \[
 Q= Q_c - \e {\mathcal{A}}^0 - \e ^2 {\mathcal{A}}^1, \quad P=P_c. 
 \]
 Dropping the subscripts $c$,   the
 bi-characteristic equations reduce to
 \begin{align*}
 \dot Q &= \nabla_P E^0 - \e \dot P \times \nabla_P \times {\mathcal{A}}' +\e^2 \nabla_P E_w, \\
 \dot P &= - \nabla_Q U - \e^2 \nabla_Q E_w.
 \end{align*}
 
 At this stage, the characteristic equations have successfully
 captured the correction term the Berry curvature. However, the
 expression for the extra wave packet energy $E_w$ is still
 complicated and may not be gauge invariant {(in fact the expression of $E_w$ is not derived in previous works \cites{GYN, GYN2})}.

{We now show that it is possible to avoid the extra
   wave packet energy if we consider well prepared initial
   condition. The expression of the extra wave packet energy is
   explicitly given by
   \begin{equation}
     \begin{aligned}
	E_w&=  \frac 12\langle \nabla_p  \chi^0, \nabla_x^2 U \nabla_p \chi^0 \rangle - \frac{\Delta_x A^0}{2A^0} + \frac{1}{2} |\nabla_x^2 S^0 \nabla_p \chi^0|^2  \\
	& \quad + \nabla_x \log A^0 \cdot {\rm Im} \left( \left\langle \chi^0, \left(p -i \nabla_z - \nabla_p E^0\right)v \right\rangle \right) - \frac{1}{2}\nabla_x^2 S^0 \nabla_p \cdot \nabla_p E^0 {\rm Im} \left( \langle \chi^0,v \rangle \right)  \\
	& \quad -\nabla_x^2 S^0 \nabla_p E^0 \cdot {\rm Im} \left( \langle \chi^0, \nabla_p v \rangle \right) + \frac{1}{2} {\rm Im} \left( \left\langle  \chi^0, (\Delta_x S^0 +2 (p-i\nabla_z)\cdot \nabla_x^2 S^0 \nabla_p)v \right \rangle \right) \\
	& \quad + {\rm Im} \left( \left\langle \chi^0, \left(p -i \nabla_z \right)\cdot \nabla_x v \right\rangle \right)+ {\rm Im} \left( \langle \chi^0,\partial_t v  \rangle \right). 
   \end{aligned}
   \end{equation}
Let us consider the case where initial condition is a plane wave
multiplied by a periodic Bloch wave,}
 \[
 A^0(0,x)=A_0,\quad S^0(0,x)=S_0 + K_0 x.
 \]
 If we denote the potential function $U(x)=c_0+c_1 x$, then, the exact solutions to equation \eqref{eqS0} and equation \eqref{eqA0} are
 \[
 A^0(t,x)=A_0.\quad S^0(t,x)=b_0(t)+b_1(t)x,
 \]
 where
 \[
 b_1(t)=K_0-c_1 t,\quad b_0(t)=S_0 - \int_0^t E^0(K_0-c_1 s) ds -c_0 t.
 \]
 In other words, $A^0(t,x)$ stays as an constant, and $S^0(t,x)$ remains be be a linear function in $x$. 
 In this case, by direct calculations,  the term $v$ reduces to $0$, and the extra wave packet energy $E_w$ 
simplifies to $0$. Hence the bi-characteristic equations become
 \begin{equation}
 \begin{cases}
 \dot Q = \nabla_P E^0 - \e \dot P \times \nabla_P \times {\mathcal{A}_{(1)}}, \\
 \dot P = - \nabla_Q U.
 \end{cases}
 \end{equation}

 This form essentially agrees with the recent results \cites{GYN,
   GYN2}, although in the above special case, the extra wave packet
 energy $E_w$ has simplified to $0$. However, our derivation is valid
 also for more general assumptions on the amplitude function and the
 phase function, and accordingly we would expect possibly different
 expression of the extra wave packet energy.  Our results show that
 the second order correction to the Bloch dynamics can be derived
 rigorously in more general situations than that considered in
 \cites{GYN, GYN2}. 
}
 
{
 \subsubsection{Linearly varying electric field}  \label{sec:quadU}

 Next, we consider the case when the potential function $U$ is
 quadratic in $x$, which corresponds to a linearly varying electric
 field. {We prove in the following that the
   characteristic equations are gauge independent if we only keep the
   static Bloch energy $\tilde{E}_s^2$ in the second order
   correction.}
 
{The bi-characteristic equation which only include
   the $\tilde{E}_s^2$ correction reduces to}
 	\begin{align}\label{eq4.2}
 	\dot{Q}&=\nabla_P E^0+\varepsilon\nabla_P(\mathcal{A}^0\cdot\nabla_Q U)+\varepsilon^2\nabla_P \tilde{E}_s^2, \\
     \label{eq4.3}
 	\dot{P}&=-\nabla_Q U-\varepsilon\nabla_Q^2 U\cdot \mathcal{A}^0-\varepsilon^2\nabla_Q \tilde{E}_s^2.
 	\end{align}

 	Notice that $\nabla_Q^2 U={\rm constant}$ when $Q$ is quadratic, and thus $\nabla_Q(\mathcal B:\nabla_Q^2 U)=  0$, then the characteristic equations become
	\begin{align*}
 	\dot{Q}&=\nabla_P E^0+\varepsilon\nabla_P \mathcal{A}^0\cdot\nabla_Q U+\varepsilon^2\nabla_P \mathcal{A}^1\cdot\nabla_Q U+\varepsilon^2\nabla_P \mathcal B:\nabla_Q^2 U, \\
 	\dot{P}&=-\nabla_Q U-\varepsilon\nabla_Q^2 U\cdot \mathcal{A}^0-\varepsilon^2\nabla_Q (\mathcal{A}^1\cdot \nabla_Q U).
 	\end{align*} 

 	We carry out the same "positional-shift" change of variable as in the last case
\begin{equation} \label{rel:change}
Q_c = Q+\varepsilon \mathcal{A}^0+\varepsilon^2 \mathcal{A}^1, \quad P_c=P, 
\end{equation}
then after length calculations (which we provide in Appendix \ref{app:b}), dropping the subscript $c$ and ignoring higher order terms, we obtain the characteristic equations in physical variables
\begin{align}
\dot{P}& =-\nabla_Q U(Q)-\varepsilon^2\nabla_Q \mathcal{A}^1\cdot\nabla_Q U,  \label{eq:phyP2} \\
\dot{Q}&= \nabla_P E^0+ \varepsilon \dot{P} \times \nabla_P \times \mathcal{A}_{(1)} - \varepsilon^2\nabla_P \mathcal{A}^0\cdot(\nabla_Q^2U\cdot \mathcal{A}^0) +\varepsilon^2 \nabla_P \mathcal B:\nabla_Q^2 U + \varepsilon^2(\nabla_Q \mathcal{A}^1)^T \cdot \nabla_P E^0.  \label{eq:phyQ2}
\end{align}
Compared with the uniform electric field case, more terms appear in the second order Bloch dynamics {due to the extra correction term of the static Bloch energy part  $\tilde{E}_s^2$} when the potential function is quadratic in $x$. 

{In fact, as shown in Appendix \ref{app:c}, all the quantities in
\eqref{eq:phyP2} and \eqref{eq:phyQ2} are gauge invariant. This suggests that when the external electric field is linearly varying,
the positional-shift change of variable can also transform the
characteristic equations in canonical variable to those in physical
variables. We remark that this scenario has not been
previously considered in the literature, and 
investigation of physical interpretations of the extra terms is interesting research directions to be pursued.}

{ 

  Finally let us make some comments on the term $E_w$ in this case, in
  general it seems challenging to prepare initial condition to make it
  vanish (as in the uniform electric field case); it is also difficult
  to establish its gauge invariance. In fact, since no a priori
  knowledge of the leading order Bloch energy $E^0$ is assumed, the
  Hamilton-Jacobi equation \eqref{eqS0} does not admit a close-form
  solution unless the potential function $U(x)$ is linear, which makes
  the analysis concerning $E_w$ intractable.  While such extra wave
  packet energy is completely ignored in previous works, we see
  through the rigorous derivation that it naturally appears in the WKB
  analysis, and thus further understanding of the term is required to
  clarify its physical implications. This is beyond the scope of our
  current work.

 }

\appendix

{
	\begin{appendices}
 {\section{Preliminaries on the Schr\"odinger equation}\label{app:a}
We follow the nondimensionalization steps as in \cite{ELY} here, as we consider the same dimensionless form \eqref{SE1}.
 
 Given a crystal lattice $\mathbb{L}$, we consider a rescaled Schr\"odinger equation \eqref{SE1}. we rewrite it below
 \begin{equation}\label{SE2}
 	i\eps \frac{ \partial}{\partial t}\psi(t,x)=H \psi(t,x)= \left(-\frac{\eps^2}{2}\Delta_x+V\left( \frac{x}{\eps}\right)+U(x)\right)\psi(t,x),
 \end{equation}
 where $V(\cdot)$ is a lattice potential, which is periodic with respect to $\mathbb{L}$ and $U(x)$ is the slow-varying scalar potential. \eqref{SE2} is a standard model for describing the motion of electrons in a perfect crystal when an external macroscopic potential is
applied. In physical units, the equation is given by
\begin{equation}\label{SE3}
	i\hbar\frac{\partial}{\partial t}\psi=-\frac{\hbar^2}{2m}\Delta \psi +V(x)\psi-U(x)\psi,
\end{equation}
where $m$ is the mass and $\hbar$ is the reduced Planck constant. As in \cite{ELY}, we introduce $l$ as the lattice constant and $\tau=ml^2/\hbar$ as the small (quantum) time scale, and  $L$ and $T$ as the macroscopic length and time scales, then
\begin{equation*}
	\begin{aligned}
	V(x)=\frac{ml^2}{\tau^2}\tilde{V}\left(\frac{x}{l}\right),\qquad U (x)=\frac{mL^2}{T^2}\tilde{U}\left(\frac{x}{L}\right). 
	\end{aligned}
\end{equation*}
By defining 
\begin{equation*}
	\begin{aligned}
	\tilde{x}=\frac{x}{L}, \qquad \tilde{t}=\frac{t}{T},\qquad \e =\frac{l}{L},\qquad h=\frac{\hbar T}{mL^2},
	\end{aligned}
\end{equation*}
one obtains after dropping the tildes,
\begin{equation}
	ih\frac{\partial}{\partial t}\psi=-\frac{h^2}{2}\Delta \psi +\frac{h^2}{\e^2}V\left(\frac{x}{\e}\right)\psi-U(x)\psi
\end{equation}
This equation has two small dimensionless parameters, $\e$ and $h$. We will only consider the
distinguished limit when $h=\e$, then we get \eqref{SE2} 

 }
 \section{The change of variable of the bi-characteristic equations in Section \ref{sec:quadU} } \label{app:b}
The change of variable  \eqref{rel:change}  for the $P$ equation is straightforward:
\begin{equation}
 	\dot{P}=-\nabla_Q U(Q_c)-\varepsilon^2\nabla_Q \mathcal{A}^1\cdot\nabla_Q U+\mathcal O(\varepsilon^3).
 	\end{equation}

 	For the $Q$ equation, notice that 
 	\begin{displaymath}
 	\dot{Q_c} = \dot{Q}+\varepsilon \dot{\mathcal{A}^0}+\varepsilon^2 \dot{\mathcal{A}^1},
 	\end{displaymath}
 	and  $\dot{P}=-\nabla_Q U(Q_c)+O(\varepsilon^2)$, we have
 	\begin{displaymath}
 	\dot{\mathcal{A}^0}=-(\nabla_P \mathcal{A}^0)^T \cdot \nabla_Q U(Q_c)+\mathcal O(\varepsilon^2),
 	\end{displaymath}
 	\begin{displaymath}
 	\dot{\mathcal{A}^1}=-(\nabla_P \mathcal{A}^1)^T \cdot\nabla_Q U(Q_c)+(\nabla_Q \mathcal{A}^1)^T\cdot \dot{Q}+\mathcal O(\varepsilon^2).
 	\end{displaymath}
 	Thus we obtain
 	\begin{align*}
 	\dot{Q_c}= & \dot{Q}+\varepsilon \dot{\mathcal{A}^0}+\varepsilon^2 \dot{\mathcal{A}^1} \\
 	= &\nabla_P E^0+\varepsilon\nabla_P \mathcal{A}^0\cdot\nabla_Q U(Q)+\varepsilon^2\nabla_P \mathcal{A}^1\cdot \nabla_Q U(Q)+\varepsilon^2 \nabla_P B:\nabla_Q^2 U \\
 	&-\varepsilon(\nabla_P \mathcal{A}^0)^T\cdot \nabla_Q U(Q_c)-\varepsilon^2 (\nabla_P \mathcal{A}^1)^T\cdot \nabla_Q U(Q_c) + \varepsilon^2(\nabla_Q \mathcal{A}^1)^T \cdot \dot{Q}+\mathcal O(\varepsilon^3).
 	\end{align*}

 Observe  that
 	\begin{displaymath}
 	\nabla_Q U(Q)=\nabla_Q U(Q_c)-\varepsilon\nabla_Q^2U(Q_c)\cdot \mathcal{A}^0 + \mathcal O(\varepsilon^2), \quad \mathcal{A}^1(P,Q) =  \mathcal{A}^1(P,Q_c)+ \mathcal O(\varepsilon),
 	\end{displaymath}
 	thus we get
 	\begin{align*}
 	\dot{Q_c}= & \nabla_P E^0+\varepsilon\nabla_P \mathcal{A}^0\cdot\nabla_Q U(Q_c)-\varepsilon^2\nabla_P \mathcal{A}^0\cdot(\nabla_Q^2U(Q_c)\cdot \mathcal{A}^0)+\varepsilon^2\nabla_P \mathcal{A}^1\cdot \nabla_Q U(Q_c)+\varepsilon^2 \nabla_P \mathcal B:\nabla_Q^2 U \\
 	&-\varepsilon(\nabla_P \mathcal{A}^0)^T\cdot \nabla_Q U(Q_c)-\varepsilon^2 (\nabla_P \mathcal{A}^1)^T\cdot \nabla_Q U(Q_c) + \varepsilon^2(\nabla_Q \mathcal{A}^1)^T \cdot \dot{Q_c},
 	\end{align*}
 	Finally, with $\dot{Q_c}= \nabla_P E^0+\mathcal O( \varepsilon)$ and some simplification, we get
 	\begin{equation}
 	\dot{Q_c}= \nabla_P E^0+ \varepsilon \dot{P} \times \nabla_P \times \mathcal{A}_{(1)} - \varepsilon^2\nabla_P \mathcal{A}^0\cdot(\nabla_Q^2U\cdot \mathcal{A}^0) +\varepsilon^2 \nabla_P B:\nabla_Q^2 U + \varepsilon^2(\nabla_Q \mathcal{A}^1)^T \cdot \nabla_P E^0+ \mathcal O(\varepsilon^3).
 	\end{equation}

 \section{Gauge invariance of some important terms}  \label{app:c}

\subsection{$\mathcal{A}^1$ is gauge invariant}

It suffices to show if we  change $\chi$ into $e^{if(p)}\chi$ and then $\mathcal{A}^1$ is still the same.

 	Recall that 
 	\begin{displaymath}
 	i \langle \chi^0,\nabla_p \chi^0\rangle \cdot \nabla_x U \chi^0 -i\nabla_p \chi^0\cdot \nabla_x U=(H^0-E^0)w,
 	\end{displaymath}
 	when we change $\chi$ into $e^{if(p)}\chi$, the left side of the equation becomes 
 	\begin{align*}
 	&i \langle e^{if(p)}\chi^0,\nabla_p e^{if(p)}\chi^0\rangle \cdot \nabla_x U e^{if(p)}\chi^0 -i\nabla_p e^{if(p)}\chi^0\cdot \nabla_x U\\
 	&\qquad=(i \langle \chi^0,\nabla_p \chi^0\rangle-\nabla_p f(p)) \cdot \nabla_x U e^{if(p)}\chi^0 -(ie^{if(p)}\nabla_p \chi^0-e^{if(p)}\chi^0\nabla_p f(p)) \cdot \nabla_x U\\
 	&\qquad=i \langle \chi^0,\nabla_p \chi^0\rangle \cdot \nabla_x U e^{if(p)}\chi^0 -ie^{if(p)}\nabla_p \chi^0 \cdot \nabla_x U\\
 	&\qquad=e^{if(p)}(i \langle \chi^0,\nabla_p \chi^0\rangle \cdot \nabla_x U \chi^0 -i\nabla_p \chi^0 \cdot \nabla_x U)\\
 	&\qquad=(H^0-E^0)e^{if(p)}w.
 	\end{align*}

 	Thus when we change $\chi$ into $e^{if(p)}\chi$, we also need to change $w$ into $e^{if(p)}w$, and $\mathcal{A}^1$ becomes the following function
 	\begin{displaymath}
 	\frac{1}{2}(i<e^{if(p)}\chi^0,\nabla_pe^{if(p)}w>+c.c.)
 	=\frac{1}{2}(i<e^{if(p)}\chi^0,e^{if(p)}\nabla_pw+if'(p)e^{if(p)}w>+c.c.).
 	\end{displaymath}
Becuase $<\chi^0,w>=0$, we conclude
 	\begin{displaymath}
 	\frac{1}{2}(i<e^{if(p)}\chi^0,\nabla_pe^{if(p)}w>+c.c.)
 	=\frac{1}{2}(i<\chi^0,\nabla_pw>+c.c.).
 	\end{displaymath}
 	
 	Hence, $\mathcal{A}^1$ is gauge independent.

\subsection{other $\mathcal O(\varepsilon^2)$ terms in the $Q$ equation are gauge invariant}
		We replace $\chi^0$ with $e^{if(P)}\chi^0$, then $\mathcal{A}^0$ becomes $\mathcal{A}^0-\nabla_P f(P)$, and $\nabla_P \mathcal{A}^0\cdot(\nabla_Q^2U \mathcal{A}^0)$ becomes	
		\[
		\nabla_P \mathcal{A}^0\cdot(\nabla_Q^2U \mathcal{A}^0)-\nabla_P \mathcal{A}^0\cdot(\nabla_Q^2U \nabla_P f(P))
		-\nabla_P^2 f(P)(\nabla_Q^2U \mathcal{A}^0)+\nabla_P^2 f(P)(\nabla_Q^2U \nabla_P f(P)).
		\]

		For convenience, we introduce the following notation,
		\[
		{\rm dif}K=K(\chi^0)-K(e^{if(P)}\chi^0),
		\]
		then we get
		\begin{align*}
		{\rm dif}\left(\nabla_P \mathcal{A}^0\cdot(\nabla_Q^2U \mathcal{A}^0)\right)=&-\nabla_P \mathcal{A}^0\cdot(\nabla_Q^2U \nabla_P f(P))-\nabla_P^2 f(P)(\nabla_Q^2U \mathcal{A}^0)+\nabla_P^2 f(P)(\nabla_Q^2U \nabla_P f(P)),
		\end{align*}
		and $\mathcal B$ becomes
		\begin{displaymath}
		-\frac{1}{4}\left(<\chi^0,\nabla_p^2\chi^0>+i\nabla_P^2f(P)-(\nabla_P f(P))^2
		+i\nabla_P f(P)<\chi^0,\nabla_P\chi^0>+i<\chi^0,\nabla_P \chi^0>^T(\nabla_P f(P))^T+c.c.\right),
		\end{displaymath}
		so we have
		\begin{displaymath}
		{\rm dif}(\mathcal B)=-\frac{1}{2}\left(-(\nabla_P f(P))^2
		+\nabla_P f(P)\mathcal{A}^0+(\mathcal{A}^0)^T(\nabla_P f(P))^T\right).
		\end{displaymath}

		In the next, we compare ${\rm dif}(\nabla_P \mathcal{A}^0\cdot(\nabla_Q^2U \mathcal{A}^0))$ with $\nabla_P\bigl({\rm dif}(\mathcal B):\nabla_Q^2 U\bigr)$. 
		In $\nabla_P(\frac{1}{2}(\nabla_P f(P))^2:\nabla_Q^2 U)$,
		we denote  $\nabla_Q^2 U=(u_{ij}),\,\nabla_P f(P)=(f_i),\,\nabla_P^2 f(P) = (f_{ij})$, then 
		\begin{displaymath}
		\frac{1}{2}(\nabla_P f(P))^2:\nabla_Q^2 U=\frac12(f_1f_1u_{11}+2f_1f_2u_{12}+2f_1f_3u_{13}+f_2f_2u_{22}+2f_2f_3u_{23}+f_3f_3u_{33}).
		\end{displaymath}
		
		Differentiating this equation with respect to $P_k$, we get
		\begin{equation}\label{A1}
		\begin{aligned}
		\partial_{P_k}(\frac{1}{2}(\nabla_P f(P))^2:\nabla_Q^2 U)=&f_{1k}f_1u_{11}+f_{1k}f_2u_{12}+f_1f_{2k}u_{12}+f_{1k}f_3u_{13}+f_1f_{3k}u_{13}\\
		&+f_{2k}f_2u_{22}+f_{2k}f_3u_{23}+f_2f_{3k}u_{23}+f_{3k}f_3u_{33} .
		\end{aligned}
		\end{equation}
		
		To compare this with $-\nabla_P^2f(p)(\nabla_Q^2U\nabla_Pf(P))$, we observe that the $k-$th component of the latter is 
		\begin{displaymath}
		-\begin{pmatrix} f_{1k} & f_{2k} & f_{3k} \end{pmatrix} \begin{pmatrix} u_{11} & u_{12} & u_{13} \\ u_{12} & u_{22} & u_{23} \\ u_{13} & u_{23} & u_{33} \end{pmatrix} \begin{pmatrix} f_{1} \\ f_{2} \\ f_{3} \end{pmatrix},
		\end{displaymath}
		which is the opposite of equation \eqref{A1}.
		
		Thus we have
		\begin{equation}\label{eqA2}
		\nabla_P(\frac{1}{2}(\nabla_P f(P))^2:\nabla_Q^2 U)-\nabla_P^2f(p)(\nabla_Q^2U\nabla_Pf(P))=0.
		\end{equation}
		
		Next, we consider $\nabla_P\bigl(-\frac12(\nabla_P f(P)\mathcal{A}^0+(\mathcal{A}^0)^T(\nabla_P f(P))^T):\nabla_Q^2 U\bigr)$. Notice that $ (\nabla_P f(P)\mathcal{A}^0)^T=(\mathcal{A}^0)^T(\nabla_P f(P))^T$ and $\nabla_Q^2 U$ is a symmetric matrix, so we have
		\begin{displaymath}
		\nabla_P f(P)\mathcal{A}^0:\nabla_Q^2 U= (\mathcal{A}^0)^T(\nabla_P f(P))^T:\nabla_Q^2 U,
		\end{displaymath}	
		\begin{displaymath}
		\begin{aligned}
		-\frac12(\nabla_P f(P)\mathcal{A}^0+(\mathcal{A}^0)^T(\nabla_P f(P))^T):\nabla_Q^2 U&=-\nabla_P f(P)\mathcal{A}^0:\nabla_Q^2 U\\
		&=-(f_1\mathcal{A}_1u_{11}+(f_1\mathcal{A}_2+f_2\mathcal{A}_1)u_{12}+(f_1\mathcal{A}_3+f_3\mathcal{A}_1)u_{13}\\
		&\qquad+f_2\mathcal{A}_2u_{22}+(f_2\mathcal{A}_3+f_3\mathcal{A}_2)u_{23}+f_3\mathcal{A}_3u_{33}).
		\end{aligned}
		\end{displaymath}
		
		Differentiating the last equation above with respect to $P_k$, we get		
		\begin{equation}
		\begin{aligned}
		\partial_{P_k}(-\nabla_P f(P)\mathcal{A}^0:\nabla_Q^2 U)=&\underbrace{-(f_{1k}\mathcal{A}_1u_{11}+(f_{1k}\mathcal{A}_2+f_{2k}\mathcal{A}_1)u_{12}+(f_{1k}\mathcal{A}_3+f_{3k}\mathcal{A}_1)u_{13}}_{B1}\\
		&\underbrace{+f_{2k}\mathcal{A}_2u_{22}+(f_{2k}\mathcal{A}_3+f_{3k}\mathcal{A}_2)u_{23}+f_{3k}\mathcal{A}_3u_{33}) }_{B1}\\
		&\underbrace{-(f_1\mathcal{A}_{1,k}u_{11}+(f_1\mathcal{A}_{2,k}+f_2\mathcal{A}_{1,k})u_{12}+(f_1\mathcal{A}_{3,k}+f_3\mathcal{A}_{1,k})u_{13}}_{B2}\\
		&\underbrace{+f_2\mathcal{A}_{2,k}u_{22}+(f_2\mathcal{A}_{3,k}+f_3\mathcal{A}_{2,k})u_{23}+f_3\mathcal{A}_{3,k}u_{33})}_{B2}.\\
		\end{aligned}
		\end{equation}	
		Here, we have needed  the notation $\mathcal{A}_{i,k}=\partial_{P_k}\mathcal{A}_i$. 

 We then compare the expression above with $\nabla_P^2 f(P)(\nabla_Q^2U \mathcal{A}^0)$ and $\nabla_P \mathcal{A}^0\cdot(\nabla_Q^2U \nabla_P f(P))$. The k-th component of $\nabla_P^2 f(P)(\nabla_Q^2U \mathcal{A}^0)$ is		
		\begin{displaymath}
		\begin{pmatrix} f_{1k} & f_{2k} & f_{3k} \end{pmatrix} \begin{pmatrix} u_{11} & u_{12} & u_{13} \\ u_{12} & u_{22} & u_{23} \\ u_{13} & u_{23} & u_{33} \end{pmatrix} \begin{pmatrix} \mathcal{A}_{1} \\ \mathcal{A}_{2} \\ \mathcal{A}_{3} \end{pmatrix},
		\end{displaymath}
		which is the opposite of part (B1).
		
		And the k-th component of $\nabla_P \mathcal{A}^0\cdot(\nabla_Q^2U \nabla_P f(P))$ is		
		\begin{displaymath}
		\begin{pmatrix} \mathcal{A}_{1,k} & \mathcal{A}_{2,k} & \mathcal{A}_{3,k} \end{pmatrix} \begin{pmatrix} u_{11} & u_{12} & u_{13} \\ u_{12} & u_{22} & u_{23} \\ u_{13} & u_{23} & u_{33} \end{pmatrix} \begin{pmatrix} f_{1} \\ f_{2} \\ f_{3} \end{pmatrix},
		\end{displaymath}
		which is the opposite of part (B2).
		
		Thus we have		
		\begin{equation}\label{A4}
		\nabla_P(-\frac12(\nabla_P f(P)\mathcal{A}^0+(\mathcal{A}^0)^T(\nabla_P f(P))^T):\nabla_Q^2 U)+\nabla_P^2 f(P)(\nabla_Q^2U \mathcal{A}^0)+\nabla_P \mathcal{A}^0\cdot(\nabla_Q^2U \nabla_P f(P))=0.
		\end{equation}
		
		From equation \eqref{eqA2} and equation \eqref{A4} we can get		
		\begin{equation}
		\nabla_P({\rm dif}(\mathcal B):\nabla_Q^2 U)-{\rm dif}(\nabla_P \mathcal{A}^0\cdot(\nabla_Q^2U \mathcal{A}^0))=0,
		\end{equation}
		which means $- \nabla_P \mathcal{A}^0\cdot(\nabla_Q^2U\mathcal{A}^0) + \nabla_P B:\nabla_Q^2 U$ is gauge independent.\\

	\end{appendices}
}


\bibliographystyle{amsxport}
\bibliography{berry,jl,newphy}

\end{document}